\documentclass[fontsize=12pt,a4paper,headings=normal,
twoside=false,leqno,parskip=half-,abstract=true]{scrartcl}
\usepackage[english]{babel}
\usepackage[utf8]{inputenc}
\usepackage{hyperref}
\hypersetup{
 bookmarks=true,
 pdftitle={meanders with three noses},
 pdfauthor={Bernold Fiedler},
 colorlinks=true,
 linkcolor=blue,
 citecolor=blue,
 filecolor=blue,
 urlcolor=blue}
  
 \usepackage{afterpage}

\usepackage{graphicx}
\usepackage[format=plain,labelfont=bf,font=small]{caption}
\usepackage{xcolor}
\usepackage[arrow, matrix, curve]{xy}
\usepackage{float}

\usepackage{caption}
\usepackage{tabulary}
\usepackage{array}
\newcolumntype{N}[1]{>{\centering\arraybackslash}m{#1}}

\usepackage{amsmath,amsthm}
\swapnumbers 
\usepackage{amssymb} 
\newcommand{\transv}{\mathrel{\text{\tpitchfork}}}
\makeatletter
\newcommand{\tpitchfork}{%
  \vbox{
    \baselineskip\z@skip
    \lineskip-.52ex
    \lineskiplimit\maxdimen
    \m@th
    \ialign{##\crcr\hidewidth\smash{$-$}\hidewidth\crcr$\pitchfork$\crcr}
  }%
}
\makeatother
\usepackage{latexsym}
\usepackage{enumerate}

\usepackage[notref,notcite,color,final 
]{showkeys}

\definecolor{refkey}{rgb}{1,0,0}
\definecolor{labelkey}{rgb}{1,0,0}

\usepackage{tikz}

  \mathchardef\ordinarycolon\mathcode`\:
  \mathcode`\:=\string"8000
  \begingroup \catcode`\:=\active
    \gdef:{\mathrel{\mathop\ordinarycolon}}
  \endgroup

\theoremstyle{plain}
\newtheorem{thm}{Theorem}[section]
\newtheorem{lem}[thm]{Lemma}
\newtheorem{prop}[thm]{Proposition}
\newtheorem{cor}[thm]{Corollary}
\newtheorem{defi}[thm]{Definition}

\begin{document}

\title{\LARGE{Design of Sturm global attractors 1:\\
Meanders with three noses, and reversibility}
\vspace{1cm}}
{\subtitle{	
	\vspace{1ex}
	{}}\vspace{1ex}
	}

\author{
 \\
\\
Bernold Fiedler* and Carlos Rocha{**}\\
\vspace{2cm}}

\date{\small{version of \today}}
\maketitle
\thispagestyle{empty}

\vfill

*\\
Corresponding author:\\
Institut für Mathematik\\
Freie Universität Berlin\\
Arnimallee 3\\ 
14195 Berlin, Germany\\
\texttt{fiedler (at) mi.fu-berlin.de}\\
 \\
{**}\\
Instituto Superior T\'ecnico\\
Avenida Rovisco Pais\\ 
1049-001 Lisboa, Portugal


\newpage
\pagestyle{plain}
\pagenumbering{roman}
\setcounter{page}{1}

\begin{abstract}
\noindent
We systematically explore a simple class of global attractors, called \emph{Sturm} due to nodal properties, for the semilinear scalar parabolic PDE
	\begin{equation*}\label{eq:*}
	u_t = u_{xx} + f(x,u,u_x) 
	\end{equation*}
on the unit interval $0 < x<1$, under Neumann boundary conditions.
This models the interplay of reaction, advection, and diffusion.

\smallskip\noindent
Our classification is based on the Sturm meanders, which arise from a shooting approach to the ODE boundary value problem of equilibrium solutions $u_t=0$.
Specifically, we address meanders with only three ``noses'', each of which is innermost to a nested family of upper or lower meander arcs.
The Chafee-Infante paradigm, with cubic nonlinearity $f=f(u)$, features just two noses.

\smallskip\noindent
Our results on the gradient-like global PDE dynamics include a precise description of the connection graphs.
The edges denote PDE heteroclinic orbits $v_1 \leadsto v_2$ between equilibrium vertices $v_1, v_2$ of adjacent Morse index.
The global attractor turns out to be a ball of dimension $d$, given as the closure of the unstable manifold $W^u(\mathcal{O})$ of the unique equilibrium with maximal Morse index $d$.
Surprisingly, for parabolic PDEs based on irreversible diffusion, the connection graph indicates time reversibility on the ($d$-1)-sphere boundary of the global attractor.

\end{abstract}

\vspace{2cm}


\textbf{
The global dynamics of the nonlinear interplay among diffusion, reaction, and advection is little understood.
This holds true even for a single equation on finite intervals, where a decreasing energy functional and nonlinear nodal properties of Sturm type considerably simplify the dynamics.
Part of this predicament is caused by an undue focus on the particular dynamics of particular nonlinearities: spatial chaos, for example, may lead to large numbers of globally competing stable and unstable equilibria.
Instead, we explore a rich class of nonlinearities with prescribed meandric equilibrium configurations of 3-nose type.
The global attractors, in that case, turn out to be balls with an attracting boundary sphere of potentially arbitrarily large dimension.
For the first time, in that class, we provide a detailed dynamic description via the global graph structure of heteroclinic orbits between equilibria.
Much to our surprise, we encountered signs of \emph{time reversibility} within the attracting boundary sphere.
This contradicts common ``knowledge'' of diffusion as \emph{the} paradigm of irreversibility.
}

\newpage
\pagenumbering{arabic}
\setcounter{page}{1}

\section{Introduction} \label{intro}

\numberwithin{equation}{section}
\numberwithin{table}{section}

The chaotic intricacies of nonautonomous second order ODE flows have been studied for many decades.
To include forced pendula, Duffing, van der Pol, Li\"enard type equations, and many others, we consider the general form
	\begin{equation}
	0 = v_{xx} + f(x,v,v_x)\,.
	\label{ODE}
	\end{equation} 
Subscripts $x$ denote derivatives of $v(x)$.

We reserve time $t$ to denote the PDE semiflow of the associated scalar reaction-advection-diffusion equation
	\begin{equation}
	u_t = u_{xx} + f(x,u,u_x)\,.
	\label{PDE}
	\end{equation}
To be specific, we consider solutions $u=u(t,x) \in \mathbb{R}$ on the unit interval $0<x<1$, with Neumann conditions $u_x=0$ at the boundaries $x=0,1$.
Subscripts $t,x$ indicate partial derivatives.
Equilibria of \eqref{PDE}, i.e. time-independent solutions $u(t,x)=v(x)$, equivalently satisfy the ``pendulum'' equation \eqref{ODE}, albeit as a boundary value problem with Neumann conditions in the spatial variable $x$.
Any chaos in \eqref{ODE} becomes spatial, in \eqref{PDE}, visible more and more prominently on longer and longer (normalized) $x$-intervals.

Many applications lead to equations of the form \eqref{PDE}, under various boundary conditions on $x$, or on the unbounded real axis.
We mention a few, cursorily.
Famous early examples include the quadratic \emph{Fisher equation} $f=\lambda u(1-u)$ of genetic selection \cite{fish}, 
and the slightly more general \emph{Kolmogorov–Petrovsky–Piskunov} (KPP) variant of population growth \cite{kpp}.
See also the stochastic branching processes addressed in \cite{bra83}.
Cubic $f$ arise in the Allen-Cahn description of interface motion in binary alloys \cite{alca79}
and, as a singular limit, in the Nagumo equation of nerve conduction.
The famous Chafee-Infante cubic $f=\lambda^2 u(1-u^2)$ falls into that class, and inspired much of the PDE analysis in the area \cite{chin74}; see also \eqref{cubic} and section \ref{ChIn} below.
The prefactor $\lambda^2$, which we now omit, arises from scaling a spatial interval $0<x<\lambda$ to unit length.
The \emph{Zeldovich–Frank-Kamenetskii equation} (ZFK) $f=u(1-u) \exp(-\beta (1-u))$ models combustion and, with the proper Arrhenius exponential instead, non-isothermal catalysis.
Chemical reactions in permeable catalysis or tubular reactors provide examples, where reaction, advection, and diffusion arise under their proper name \cite{aris}.
Quasilinear variants of \eqref{PDE} arise, for example, in curve shortening and interface flows \cite{an91, figt04, figt06}.
Many applications involve singular limits.
For applications to viscous hyperbolic balance laws, see for example \cite{hae99}.
A few of the formidable complications of $x$-dependent nonlinearities $f$ have been tackled with in \cite{ampp};
see also \cite{fietal02,haewol}.
The PDE \eqref{PDE} also appears as a parabolic limit in problems of, both, elliptic and hyperbolic type
	\begin{equation}\label{PDEeh}
	\pm \varepsilon^2 u_{tt}+u_{xx} -u_t+ f(x,u,u_x) = 0\,,
	\end{equation}
when dominated by advection or damping, respectively \cite{mie94,schee96,moso89,fiscvi}.
A spatially discrete variant models a coupled chain of overdamped pendula \cite{fibh}; see also subsection \ref{jac}.
Conley index theory, as a homotopy-invariant, global topological tool, has extended the Chafee-Infante paradigm further to include applications to certain beam equations, and settings like FitzHugh-Nagumo, Cahn-Hilliard, and certain phase field equations \cite{hami91,mi95}.
See \cite{fife} for a broad earlier survey on phase field equations.
More recently, and mostly for systems of equations in biological context, see \cite{murr}.
See also the survey \cite{fisc03} for further mathematical and applied aspects.
In the spirit of \eqref{PDE}, very interesting global results for Ginzburg-Landau patterns on 2-spheres, and other compact surfaces of revolution, have recently been obtained by \cite{dai21,daila21}.
Meanwhile, the mathematical literature on reaction-diffusion equations alone, as refereed in Zentralblatt under MSC 35K57, has grown to more than 15,000 entries \cite{zb}.

It is therefore not our intention, in the present paper, to contribute just another analysis or simulation, for this or that particular nonlinearity $f$, arising in one or the other highly specialized applied context.
For general $x$-dependent nonlinearities, on the other hand, the chaotic complexities of even the ODE equilibrium problem \eqref{ODE} seem to frustrate any all-out attack on the PDE dynamics of \eqref{PDE}, a priori.
Or, do they?

In fact, it is possible to characterize the class of all ODE equilibrium ``configurations'', qualitatively, by certain permutations $\sigma$.
See the following section \ref{Back}.
The permutations $\sigma$ themselves, as introduced by Fusco and Rocha \cite{furo91}, are based on the discrepancies between the orderings of the equilibria at the boundaries $x=0$ and  $x=1$, respectively; see \eqref{hdef},\eqref{permdef} below.
Although each of the permutations will be represented by an open class of nonlinearities $f$, in principle, we will provide specific  nonlinearities only in exceptional cases; but see \eqref{cubic} and section \ref{ChIn} for cubic $f$.
In general,
\begin{quotation}
\noindent
\emph{it will therefore be the qualitative configuration of ODE equilibria \eqref{ODE}, which we assume to be given, rather than some particular nonlinearity $f$.}
\end{quotation}
Our approach is somewhat reminiscent of algebra: it is much easier to construct a polynomial with given zeros, you know, than to determine all zeros of a polynomial.

In the present paper, we describe the global dynamics of the full PDE \eqref{PDE}, for a certain subclass of permutations $\sigma$.
This allows us to design certain time asymptotic global attractors of \eqref{PDE}, with three competing attracting sinks.
A plethora of other equilibria, of arbitrarily high unstable dimension, may be involved in the boundaries of their domains of attraction.
The resulting PDE dynamics turns out to be gradient-like, by a general energy functional.
In particular, the PDE dynamics on the global attractor will consist of equilibria and their heteroclinic orbits \eqref{het}, only.
Still, we will encounter at least some of the intricacies which are caused by the competition among large numbers of highly unstable equilibria.

Some mathematical generalities are easily settled.
For continuously differentiable nonlinearities $f\in C^1$, standard theory of strongly continuous semigroups provides local solutions $u(t,x)$ of \eqref{PDE} in suitable Sobolev spaces $u(t, \cdot) \in X \subseteq C^1 ([0,1], \mathbb{R})$, for $t \geq 0$ and given initial data $u=u_0(x)$ at time $t=0$.
See \cite{he81, pa83, ta79} for a general PDE background.

We assume the solution semigroup $u(t,\cdot)$ generated by the nonlinearity $f$ to be \emph{dissipative}: any solution $u(t,\cdot)$ exists globally in forward time $t\geq 0$, and eventually enters a fixed large ball in $X$.
Explicit sufficient, but by no means necessary, conditions on $f=f(x,u,p)$ which guarantee dissipativeness are sign conditions $f(x,u,0)\cdot u<0$, for large $|u|$, together with subquadratic growth in $|p|$.
For large times $t\rightarrow\infty$, any large ball in $X$ then limits onto the same maximal compact and invariant subset $\mathcal{A}=\mathcal{A}_f$ of $X$ which is called the \emph{global attractor}. 
In general, the global attractor $\mathcal{A}$ consists of all solutions $u(t,\cdot)$ which exist globally, for all positive and negative times $t\in\mathbb{R}$, and remain bounded in $X$.
Of course, $\mathcal{A}$ therefore contains any equilibria, heteroclinic orbits, basin boundaries, or more complicated recurrence which might arise, in general.
See \cite{bavi92, chvi02, edetal94, ha88, haetal02, la91, ra02, seyo02, te88} for global attractors in general.

In the specific setting \eqref{PDE}, which possesses much additional structure, we call the global attractors $\mathcal{A}$ \emph{Sturm}.
The beautiful survey \cite{ra02} puts some previous work on Sturm attractors in a broader perspective.
See also \cite{firo3d-1, firo3d-2, firo3d-3} and the references there.
So far, for general theory.

\section{Background and outline}\label{Back}
Admittedly, the above information on Sturm attractors $\mathcal{A}$ is quite general.
However, it provides practically no information concerning the specific dynamics on $\mathcal{A}$.
Rather than complacently pontificate a few pretty vague generalities, here, we aim to elucidate at least some of that very rich inner dynamics.
Already the chaotic intricacies of the mere equilibrium ODE \eqref{ODE} may hint at the scope of our quest.
In particular, after decades of dedication and quite a few unexpected results, we hope to convince our readers that the purportedly ``trivial'' dynamics of \eqref{PDE} is still poorly understood.
That is why we proceed by examples.

\begin{figure}[p!]
\centering \includegraphics[width=\textwidth]{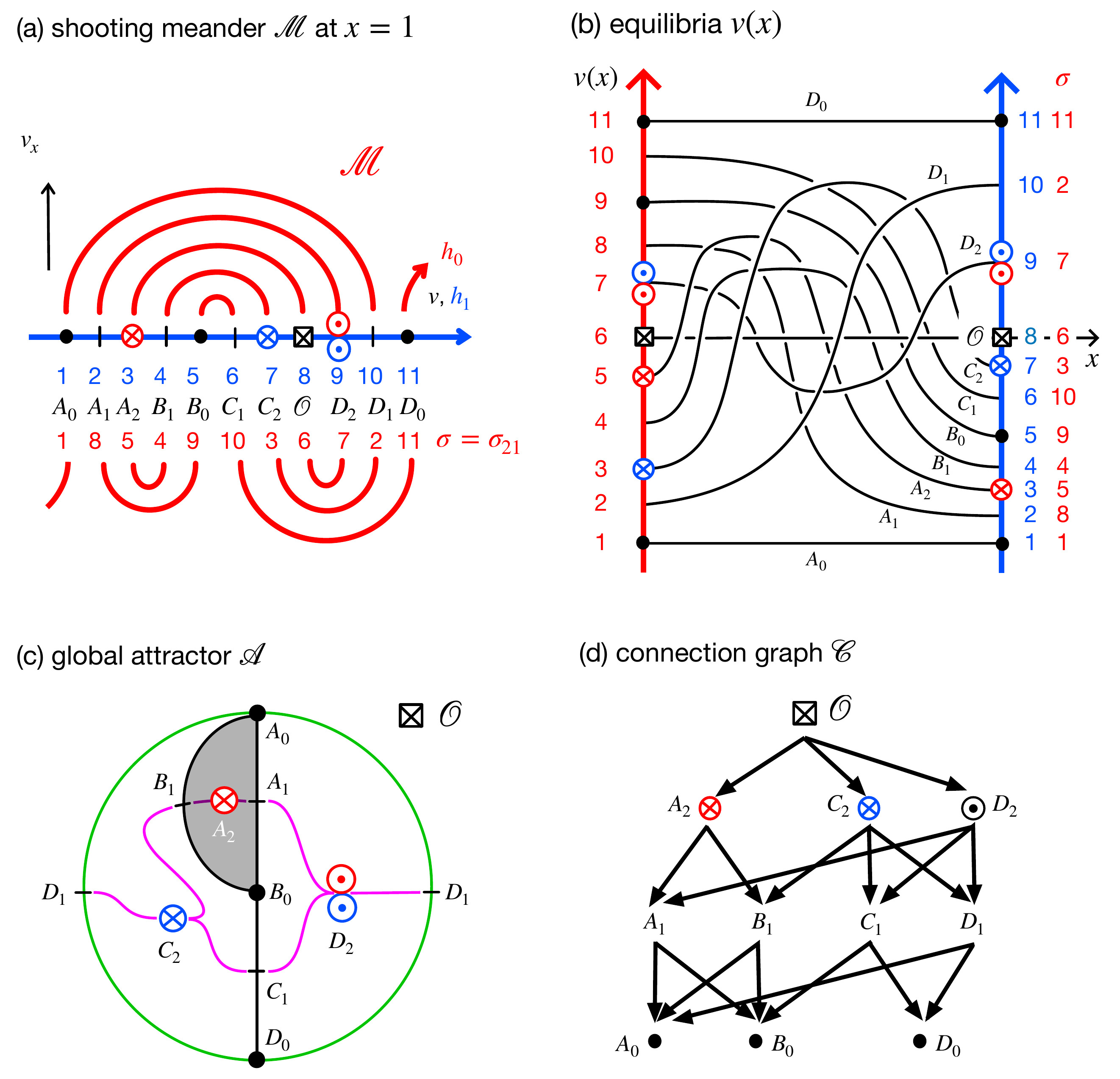}
\caption{\emph{
Example of a Sturm 3-ball global attractor $\mathcal{A}= \mathrm{clos}\ W^u(\mathcal{O})$. 
Equilibria \hbox{are labeled} as $\mathcal{E}=
\{A_0,A_1,A_2,B_1,B_0,C_1,C_2,\mathcal{O},D_2,D_1,D_0\}$. 
Subscripts and Morse indices coincide; $i(\mathcal{O})=3$.
Black dots mark sinks $i=0$, and small annotated circles (red, blue) indicate $i=2$.
The previous papers \cite{firo3d-1, firo3d-2} established the equivalence of the viewpoints (a)--(c).
(a) The Sturm meander $\mathcal{M}$ of the global attractor $\mathcal{A}$. 
The meander $\mathcal{M}$ is a stylized representation of the curve $a \mapsto (v,v_x)$, at $x=1$, which results from Neumann initial conditions $(v,v_x)=(a,0)$, at $x=0$, by shooting via the equilibrium ODE \eqref{ODE}. 
Intersections of the meander with the horizontal $v$-axis indicate equilibria. 
Transversality of intersections at $v=v(1)$ is equivalent to hyperbolicity of $v$.
(b) Spatial profiles $x\mapsto v(x)$ of the equilibria $v \in \mathcal{E}$. 
Note the different orderings of $v(x)$, by $h_0$ at the left boundary $x=0$ (red), and by the Sturm permutation $\sigma = h_0^{-1}h_1$ at the right boundary $x=1$. 
The same orderings characterize the meander $\mathcal{M}$ in (a).
(c) The Thom-Smale or Sturm complex $\mathcal{S}$ of the boundary 2-sphere $\Sigma^2=\partial \mathcal{A}=\partial W^u(\mathcal{O})$.
The right and left boundaries $A_0D_1D_0$ (green) have to be identified with each other. 
Note the 2-cells given by the unstable manifolds of the three equilibria $A_2,C_2,D_2$ of Morse index $i=2$.
The Chafee-Infante 2-cell of $A_2$ is shaded gray.
The basin boundaries of the stable equilibria $A_0,B_0,D_0$ in the 2-sphere are indicated in purple.
(d) The connection graph $\mathcal{C}$ of $\mathcal{A}$.
Vertices are the equilibria, ranked by rows of descending Morse index.
Directed edges indicate heteroclinic orbits $\leadsto$ between Morse-adjacent equilibria.
Comparison with the Thom-Smale complex (c) verifies transitivity and cascading of heteroclinic orbits.
}}
\label{3ball}
\end{figure}

Two additional structures help, in our Sturm setting.
First, \eqref{PDE} possesses a \emph{Lyapunov~function}, alias a variational or gradient-like structure, under separated boundary conditions;  see \cite{ze68, ma78, mana97, hu11, fietal14, lafi18, labe22}. 
Therefore the time invariant global attractor consists of equilibria and of solutions $u(t, \cdot )$, $t \in \mathbb{R}$, with forward and backward limits, i.e.
	\begin{equation}
	\underset{t \rightarrow -\infty}{\mathrm{lim}} u(t, \cdot ) = v_1\,,
	\qquad
	\underset{t \rightarrow +\infty}{\mathrm{lim}} u(t, \cdot ) = v_2\,.
	\label{het}
	\end{equation}
In other words, the $\alpha$- and $\omega$-limit sets of $u(t,\cdot )$ are two distinct equilibria $v_1$ and $v_2$.
We call $u(t, \cdot )$ a \emph{heteroclinic} or \emph{connecting} orbit, or \emph{instanton},  and write $v_1 \leadsto v_2$ for such heteroclinically connected equilibria. 
See fig.~\ref{3ball}(c),(d) for a modest 3-ball example with $N=11$ equilibria.
Although the variational structure persists for other separated boundary conditions, the possibility of rotating waves shows that it may fail under periodic boundary conditions.
See however \cite{fietal14,firowo12}.

The second structure is a \emph{Sturm nodal property}, which we express by the \emph{zero number} $z$.
Let $0 \leq z (\varphi) \leq \infty$ count the number of (strict) sign changes of continuous spatial profiles $\varphi : [0,1] \rightarrow \mathbb{R}, \, \varphi \not\equiv 0$.
For any two distinct solutions $u^1$, $u^2$ of \eqref{PDE}, the zero number
	\begin{equation}
	t \quad \longmapsto \quad z(u^1(t, \cdot ) - u^2(t, \cdot ))\,
	\label{zdrop}
	\end{equation}
is then nonincreasing with time $t$, for $t\geq0$, and finite for $t>0$.
Moreover $z$ drops strictly, with increasing $t>0$, at any multiple zero of the spatial profile $x \mapsto u^1(t_0 ,x) - u^2(t_0 ,x)$; see \cite{an88}.
This remains true under other separated or periodic boundary conditions.
See Sturm \cite{st1836} for the linear autonomous variant. 

The consequences of the Sturm nodal property \eqref{zdrop} for the nonlinear dynamics of \eqref{PDE} are enormous.
For an introduction see \cite{ma82, brfi88, fuol88, mp88, brfi89, ro91, fisc03, ga04} and the many references there.
Already Sturm observed that all eigenvalues $\mu_0>\mu_1>\ldots$ of the PDE linearization of \eqref{PDE} at any equilibrium $v$ are algebraically simple and real.
In fact $z(\varphi_j)=j$, for the eigenfunction $\varphi_j$ of $\mu_j$.
\emph{We assume all equilibria are hyperbolic, i.e. all eigenvalues are nonzero.}
The \emph{Morse index} $i(v)$ of $v$ then counts the number of unstable eigenvalues $\mu_j>0$.
In other words, the Morse index $i(v)$ is the dimension of the unstable manifold $W^u(v)$ of $v$.
Let $\mathcal{E} \subseteq \mathcal{A}$ denote the set of equilibria.
Our generic hyperbolicity assumption and dissipativeness of $f$ imply that $N$:= $|\mathcal{E}|$ is odd; see also \eqref{euler}.

Surprisingly, Morse-Smale transversality, a prominent concept in \cite{pasm70, pame82, ol83}, is an automatic nonlinear consequence of \eqref{zdrop}, given hyperbolicity of equilibria \cite{he85, an86}.
More precisely, intersections of unstable and stable manifolds $W^u(v_1)$ and $W^s(v_2)$ along heteroclinic orbits $v_1 \leadsto v_2$ are  automatically transverse:
$W^u(v_1) \transv W^s(v_2)$.
In the Morse-Smale setting, Henry \cite{he85} also observed
	\begin{equation}
	\label{eq:1.5b}
	v_1 \leadsto v_2 \quad\Longleftrightarrow\quad v_2\in\partial W^u(v_1)\,.
	\end{equation} 
Here $\partial W^u(v) := \mathrm{clos}\,W^u(v) \setminus W^u(v)$ denotes the topological boundary of the unstable manifold $W^u(v)$.

In a series of papers, based on the zero number, we have given a purely combinatorial description of Sturm global attractors $\mathcal{A}$; see \cite{firo96, firo99, firo00}.
Define the two \emph{boundary orders} $h_0, h_1$: $\lbrace 1, \ldots, N \rbrace \rightarrow \mathcal{E}$ of the equilibria such that
	\begin{equation}
	h_\iota (1) < h_\iota (2) < \ldots < h_\iota (N) \qquad \mathrm{at}
	\qquad x=\iota \in \{0,1\}\,.
	\label{hdef}
	\end{equation}
See fig.~\ref{3ball}(b) for an example with $N=11$ equilibrium profiles $v(x)$, enumerated by labels $\mathcal{E} = \{1,\ldots,11\},$ and boundary orders $h_0 = \mathrm{id},\ h_1 = (1\ 10\ 7\ 4\ 3\ 8\ 9\ 5\ 6\ 11)$.	
The general combinatorial description of Sturm global attractors $\mathcal{A}$ is based on the \emph{Sturm~permutation} $\sigma \in S_N$\,, defined by Fusco and Rocha in \cite{furo91} as
	\begin{equation}
	\sigma:= h_0^{-1} \circ h_1\,.
	\label{permdef}
	\end{equation}
Already in \cite{ro91}, the following explicit recursions have been derived for the Morse indices $i_j:=i(h_0(j))$ along the meander:
	\begin{equation}
	\begin{aligned}
	i_1 &=  i_N = 0\,;\\
	i_{j+1} &= i_j+(-1)^{j+1}\,
	\mathrm{sign}\, (\sigma^{-1}(j+1)-\sigma^{-1}(j))\,.\\
	\end{aligned}
	\label{i}
	\end{equation}
The zero numbers, $z_{jk} := z(h_0(j)-h_0(k))\geq 0$ for $j\neq k$, are given recursively by
	\begin{equation}
	\begin{aligned}
	z_{kk} &:= i_k\,;\qquad\qquad\qquad \\
	z_{1k} &\phantom{:}= z_{Nk}=0\,;\\
	z_{k+1,k}&\phantom{:}= \min\{i_k\,,i_{k+1}\};\\
	z_{j+1,k} &\phantom{:}= \scalebox{0.97}[1.0]
	{$ z_{jk} + \tfrac{1}{2}(-1)^{j+1}
	\cdot\big[ \mathrm{sign}\,\big(\sigma^{-1}(j+1)-\sigma^{-1}(k)\big)-\mathrm{sign}\,\big(\sigma^{-1}(j)-\sigma^{-1}(k)\big)\big].
	$}
	\end{aligned}
	\label{z}
	\end{equation}

Using a shooting approach to the ODE boundary value problem \eqref{ODE}, the Sturm~permutations $\sigma \in S_N$ have been characterized, purely combinatorially, as \emph{dissipative Morse meanders} in \cite{firo99}.
Here the \emph{dissipativeness} property, abstractly, requires fixed $\sigma(1)=1$ and $\sigma(N)=N$.
In fact, the shooting meander emanates upwards, towards $v_x>0$, from the leftmost (or lowest) equilibrium at $\sigma(1)=1$, and terminates from below, $v_x<0$, at x=1.
The \emph{meander} property requires the formal path $\mathcal{M}$ of alternating upper and lower half-circle arcs defined by the permutation $\sigma$, as in fig.~\ref{3ball}(c), to be Jordan, i.e.~non-selfintersecting.
For dissipative meanders, the recursion in \eqref{i}, and $i_1=0$, define all Morse numbers $i_j$\,.
Note how $j$ and $i_j$ are always of opposite parity, $\mod 2$.
In particular, $N$ is odd, and $i_N$=0 follows automatically.
The \emph{Morse} property, finally, requires nonnegative Morse indices $i_j\geq0$ in the formal recursion \eqref{i}, for all $j$.
For brevity, we also use the term \emph{Sturm meanders}, for dissipative Morse meanders.

For a simple recipe to determine the Morse property of a meander, \emph{the Morse number increases by 1, along any right turning meander arc, but decreases by 1 along left turns.}
This holds, independently, for upper and lower meander arcs, and remains valid even when the proper orientation of the arc is reversed; see \eqref{i}.
For examples see figs.~\ref{3ball}, \ref{susp}, \ref{ci3}, and \ref{m22susp}.
The beautifully illustrated book \cite{ka17} contains ample material on many additional aspects of meanders.
Even ``just'' counting meanders, with a prescribed number of ``noses'' \eqref{nose}, is a deep and fascinating subject \cite{del18, zog20}.
The results for Morse meanders are much less explicit, so far \cite{wo17}.

In the present paper, we address Sturm meanders.
We will return to the intriguing issue of non-Morse dissipative meanders with some negative ``Morse indices'' $i_j<0$, briefly, in proposition \ref{suspprop} and section \ref{Nonmorse}.
See also our sequel \cite{firo23}.

More geometrically, global Sturm attractors $\mathcal{A}_f$ and $\mathcal{A}_g$ of dissipative nonlinearities $f, g$ with the same Sturm permutation $\sigma_f = \sigma_g$ are $C^0$ orbit-equivalent \cite{firo00}.
Only locally, i.e.~for $C^1$-close nonlinearities $f$ and $g$, this global rigidity result is based on the Morse-Smale transversality property mentioned above. See for example \cite{pasm70, pame82, ol83}, for such local aspects.
Section \ref{Rot} discusses some ``trivial equivalences'' between Sturm attractors $\mathcal{A}_f$ and $\mathcal{A}_g$ with different Sturm permutations $\sigma_f \neq \sigma_g$\,.

In \cite{firo96} we have shown how to determine which equilibria $v_1,v_2$ possess a heteroclinic orbit connection \eqref{het}, explicitly and purely combinatorially from dissipative Morse meanders $\sigma$. 
In the elegant formulation of Wolfrum  \cite{wo02},
\begin{equation}\label{wolfrum}
v_1\leadsto v_2 \quad\Longleftrightarrow\quad  v_1, v_2\  \textrm{are}\ z\textrm{-adjacent, and}\ i(v_1)>i(v_2)\,;
\end{equation}
see also the comment in the appendix of \cite{firo3d-2}.
Here equilibria $v_1\neq v_2$ are called $z$-\emph{adjacent}, if there does not exist any \emph{blocking equilibrium} $w$ strictly between $v_1$ and $v_2$, at $x=0$ (or, equivalently, at $x=1$) such that
\begin{equation}\label{block}
z(v_1-w)=z(w-v_2)=z(v_1-v_2).
\end{equation}
With \eqref{wolfrum}, all heteroclinic orbits then follow from \eqref{i} and \eqref{z} above.

Clearly, any heteroclinic orbit $u(t,.): v_1\leadsto v_2$ implies adjacency: by \eqref{zdrop}, any blocking equilibrium $w$ would force $z(u(t,.)-w)$ to drop strictly at the Neumann boundary $x=0$, for some $t=t_0$.
This contradicts the equal values of $z$ at the limiting equilibria $v_1,v_2$ of $u$, for $t\rightarrow\pm\infty$.

As a trivial corollary, for example, we conclude $v_1\leadsto v_2$, for neighbors $v_1,v_2$ on any boundary order $h_\iota$. Here we label $v_1,v_2$ such that $i(v_1)=i(v_2)+1$; see \eqref{i}.
For an in-depth analysis and many more examples see \cite{rofi21}.

We encode the above heteroclinic structure in the directed \emph{connection graph} $\mathcal{C}$.
See fig.~\ref{3ball}(d) for an example.
The connection graph is graded by the Morse index $i$ of its equilibrium vertices.
Directed edges are the heteroclinic orbits $v_1\leadsto v_2$ running downwards between equilibria of adjacent Morse index.
Uniqueness of such heteroclinic orbits, given $v_1,v_2$, had already been observed in lemma 3.5 of \cite{brfi89}; see also \cite{furo91}.

Directed paths in the connection graph in fact encode all heteroclinic orbits.
Indeed, the heteroclinic relation $\leadsto$ on $\mathcal{E}$ is transitive, by Morse-Smale transversality and the $\lambda$-Lemma \cite{pame82}.
Therefore, any directed path from $v_1$ to $v_2$ also defines a direct heteroclinic orbit $v_1\leadsto v_2$.
Given $v_1\leadsto v_2$, conversely, the cascading principle first described in \cite{brfi89} asserts an interpolating sequence of heteroclinic orbits between equilibria of adjacent Morse indices, from $v_1$ to $v_2$\,. 

The \emph{basin of attraction} of an $i=0$ sink vertex $v$ in $\mathcal{A}$, for example, consists of $v$ itself, and all heteroclinic orbits $v_1\leadsto v$. 
The \emph{basin boundary} consists of just those other equilibria $v_1$, and all heteroclinic orbits among them.
The connection graph $\mathcal{C}$ readily identifies all those equilibria.
See our discussion in reversibility subsection \ref{Dis-rev} for a nontrivial geometric example \eqref{revrev},\eqref{revequi} based on the connection graph of fig.~\ref{conngraph}(c).

Recently, we have embarked on a more explicitly geometric description of Sturm attractors.
The disjoint dynamic decomposition
	\begin{equation}
	\mathcal{A} = \bigcup\limits_{v \in \mathcal{E}} W^u(v) 
	\label{S}
	\end{equation}
of the global attractor $\mathcal{A}$ into unstable manifolds $W^u(v)$ of equilibria $v$ is called the \emph{Thom-Smale complex} or \emph{dynamic complex}; see for example \cite{fr79, bo88, bizh92}.
In our Sturm setting \eqref{PDE} with hyperbolic equilibria $v \in \mathcal{E}$\,, the Thom-Smale complex is a finite regular cell complex, in the terminology of algebraic topology: the boundaries $\mathrm{clos}\,W^u\setminus W^u$ of the open $i(v)$-cells $W^ u (v)$ are homeomorphic to spheres of dimension $i(v)-1$.
The proof follows from the Schoenflies property \cite{firo14,firo15}.
We therefore call the regular cell decompositions \eqref{S} of the Sturm global attractor $\mathcal{A}$ the \emph{Sturm~complex} $\mathcal{S}$.

We call $d=\dim \mathcal{A} := \max_{v \in \mathcal{E}} i(v)$ the \emph{dimension} of $\mathcal{A}$, or of the complex $\mathcal{S}$.
Then at least one equilibrium $\mathcal{O}$ has maximal Morse index $i(\mathcal{O})=d$, i.e. $i(v) \leq d$ for all other Morse indices.
If $\mathcal{A} = \mathrm{clos}\,W^u(\mathcal{O})$ is the closure of a single $d$-cell, then the Sturm complex turns out to be a closed $d$-ball \cite{firo15}.
We call this case a \emph{Sturm $d$-ball}.

A 3-dimensional Sturm complex $\mathcal{C}$\,, for example, is the regular Thom-Smale complex of a 3-dimensional Sturm global attractor $\mathcal{A}$\,.
See fig.~\ref{3ball}(c) for the Sturm complex $\mathcal{S}$ of the Sturm 3-ball $\mathcal{A}$ associated to the meander in fig.~\ref{3ball}(a).

In the Sturm-ball trilogy \cite{firo3d-1, firo3d-2, firo3d-3} we have characterized all Sturm 3-balls $\mathcal{S}$.
Earlier, the trilogy \cite{firo2d-2, firo2d-1, firo2d-3} had characterized all planar Sturm complexes $\mathcal{S}$, i.e. the case $\dim \mathcal{A}=2$.
The case $\dim \mathcal{A}=1$, i.e.~$\sigma=\mathrm{id}_N$ with odd $N\geq 3$, is a trivial line with alternating $i=0$ sinks and $i=1$ saddles.
Global asymptotic stability of a unique sink equilibrium is the case $N=1$ of $\dim \mathcal{A}=0$.
 
Conversely, we have described in \cite{firo20, rofi21} how the boundary orders $h_\iota$ of \eqref{hdef}, and therefore the Sturm permutation $\sigma$ of \eqref{permdef}, are determined uniquely by the \emph{signed hemisphere decomposition}.
This is a slight refinement of the Sturm complex $\mathcal{S}$, which we do not pursue in further detail here.
In fig.~\ref{3ball}, for example, the signed hemisphere complex (c) determines how the boundary orders $h_0$ (red in (a)) and $h_1$ (blue) traverse the equilibrium vertices, from the North pole $A_0$ to the South pole $D_0$.
The predecessors and successors, on $h_\iota$\,, of the repelling sphere barycenter $\mathcal{O}$ are marked by small annotated red and blue circles, everywhere in fig.~\ref{3ball}.

The above results have illustrated the central importance of the Sturm permutations $\sigma$ or, equivalently, their Sturm meanders $\mathcal{M}$, for a systematic description of Sturm global attractors $\mathcal{A}$ and their Sturm complexes $\mathcal{S}$.
In the present paper we discuss Sturm attractors which arise from Sturm meanders with at most three \emph{noses} (called ``pimples'', in \cite{zog20}; see also the (2,1)-lieanders in \cite{del18}).
Here noses are subscripts $j\in\{1,\ldots,N-1\}$ such that
\begin{equation}
\label{nose}
\sigma(j+1)=\sigma(j)\pm 1.
\end{equation}
In other words, the associated meander vertices are adjacent under, both, $h_0$ and $h_1$.

The simplest case, of just two noses, is called the \emph{Chafee-Infante attractor}.
This has been well-studied, ever since it first arose for cubic nonlinearities $f=\lambda^2 u(1-u^2)$ in \cite{chin74}.
As a warm-up on terminology, and as a simple illustration of our approach, we review this case in section \ref{ChIn}.
For a 3-nose meander see fig.~\ref{3ball} again.

Section \ref{Res} then presents our main results on the general case of primitive 3-nose meanders $\mathcal{M}_{pq}$ with two nose arcs above the horizontal axis, each as the innermost of $p$ and $q$ nested upper arcs, respectively.
Below the horizontal axis, the only remaining nose is centered as the innermost of the complementing $p+q$ lower arcs.
Since all lower arcs are nested, we also call that configuration a (lower) \emph{rainbow}.
It turns out that the resulting curves are meanders if, and only if, $p-1$ and $q+1$ are co-prime, i.e., they do not share any nontrivial integer factor.
See theorems \ref{nonmorse}, \ref{krmuthm}, where it is also established that the dissipative meander $\mathcal{M}_{pq}$ is Sturm if, and only if, $p=r(q+1)$, for some $r,q\geq1$.
Let $\sigma_{rq}$ denote the associated Sturm permutations.
The resulting global attractors $\mathcal{A}_{rq}$ are all distinct -- except for the not immediately obvious ``trivial'' linear flow equivalence upon interchange of $r$ and $q$; see corollary \ref{krrhocor}.
In theorem \ref{krballthm}, the Sturm complex $\mathcal{S}_{rq}$ turns out to be a Sturm ball of dimension $r+q$.
The 3-ball attractor of fig.~\ref{3ball}, for example, is trivially equivalent to the simple case $r=2,\ q=1$, in the sense of section \ref{Rot}.

Quite surprisingly, the connection graph $\mathcal{C}_{rq}$\,, restricted to the invariant boundary sphere of the Sturm ball $\mathcal{A}_{rq}$\,, turns out to be time reversible; see section \ref{krrevpf}.
Although this is also true in the Chafee-Infante case, it is a quite unexpected phenomenon in parabolic diffusion equations which most of us would rightly consider \emph{the} paradigm of irreversibility.
\emph{Time reversibility} in its strongest form means the existence of an involutive \emph{reversor} $\mathcal{R}: \Sigma \rightarrow \Sigma$ which reverses the time direction of PDE orbits of \eqref{PDE}, on a ``large'' invariant subset $\Sigma\subset\mathcal{A}$.
In particular, with any two equilibria $v_1,v_2 \in \Sigma$ such that $v_1 \leadsto v_2$ in $\mathcal{A}$, the subset $\Sigma$ should also contain some of those heteroclinic orbits.
Restricted to equilibria $v_1,v_2\in \mathcal{E}\cap \Sigma$, strong reversibility implies the weaker statement
\begin{equation}
\label{rev}
v_1\leadsto v_2 \quad\Longleftrightarrow\quad \mathcal{R}v_2 \leadsto \mathcal{R}v_1
\end{equation}
on $\Sigma$. In other words, the reversor $\mathcal{R}$ induces an automorphism of the connection di-graph $\mathcal{C}$, on the vertices in $\Sigma$ and their edges, which reverses edge orientation.
The connection graph of fig.~\ref{3ball}(d), for example, illustrates reversibility \eqref{rev} under the reversor
\begin{equation}
\label{3ballrev}
\mathcal{R}: \quad A_j \longleftrightarrow D_{2-j}\,,\quad B_j \longleftrightarrow C_{2-j}\,,
\end{equation}
of the equilibria on the boundary 2-sphere $\Sigma:=\partial W^u(\mathcal{O})=\partial\mathcal{A}$ of fig.~\ref{3ball}(c).

We prove theorem \ref{nonmorse} in section \ref{Nonmorse}.
To circumvent tiresome mathematical pedantry, we only provide proofs for the simplest interesting case $r=1,\ p=q+1$ of our remaining results, in section \ref{r=1}.
This includes the explicit connection graphs $\mathcal{C}_{1q}$ for $q\geq2$; see theorem \ref{connthm}.

Section \ref{Dis} touches the general case $\sigma_{rq}$\,, which will be addressed in our sequel \cite{firo23}.
We also discuss some non-dissipative PDE aspects, and a spatially discrete ODE variant of \eqref{PDE}.
We conclude with more geometric ODE models of the connection graphs $\mathcal{C}_{1q}$ and their time reversibility.

\section{Rotations, inverses, and suspensions}\label{Rot}
To reduce the sheer number of cases, a proper consideration of symmetries is mandatory.
In this section we recall the notion of trivial equivalence for Sturm attractors $\mathcal{A}$, meanders $\mathcal{M}$, permutations $\sigma$, and connection graphs $\mathcal{C}$, as introduced in \cite{firo96}; see also section 3 in \cite{firo3d-3}.
As a prelude to induction over the number of arcs in 3-nose meanders, we also discuss double cone suspensions $\widetilde{\mathcal{A}},\widetilde{\mathcal{M}},\widetilde{\sigma},\widetilde{\mathcal{C}}$ of the entourage $\mathcal{A},\mathcal{M},\sigma,\mathcal{C}$. See also previous accounts in \cite{firo00,ka17,rofi21}.

\emph{Trivial equivalences} are generated as the Klein 4-group $\langle \kappa,\rho \rangle$ with commuting involutive generators
	\begin{align}
	(\kappa u)(x) &:= -u(x)\,; \label{rot}\\
	(\rho u)(x) &:= u(1-x)\,.\label{inv}
	\end{align}
In the PDE \eqref{PDE}, the $u$-flip $\kappa$ induces a linear flow equivalence $\kappa: \mathcal{A}_f \rightarrow \kappa \mathcal{A}_f = \mathcal{A}_{f^\kappa}$ of the global attractors with nonlinearities $f(x,u,p)$ and  $f^\kappa (x,u,p)$:= $f(x,-u,-p)$.
Similarly, $x$-reversal $\rho$ induces a linear flow equivalence $\rho: \mathcal{A}_f \rightarrow \rho \mathcal{A}_f = \mathcal{A}_{f^\rho}$ via $f^\rho (x,u,p)$:= $f(1-x,u,-p)$. Here and below $\mathcal{A},\ \mathcal{E},\ \mathcal{M},\ h_\iota,\ \sigma$ refer to $f$, whereas $\mathcal{A}^\gamma,\ \mathcal{E}^\gamma,\ \mathcal{M}^\gamma,\ h_\iota^\gamma,\ \sigma^\gamma$ will refer to $f^\gamma$.

\begin{figure}[t!]
\centering \includegraphics[width=\textwidth]{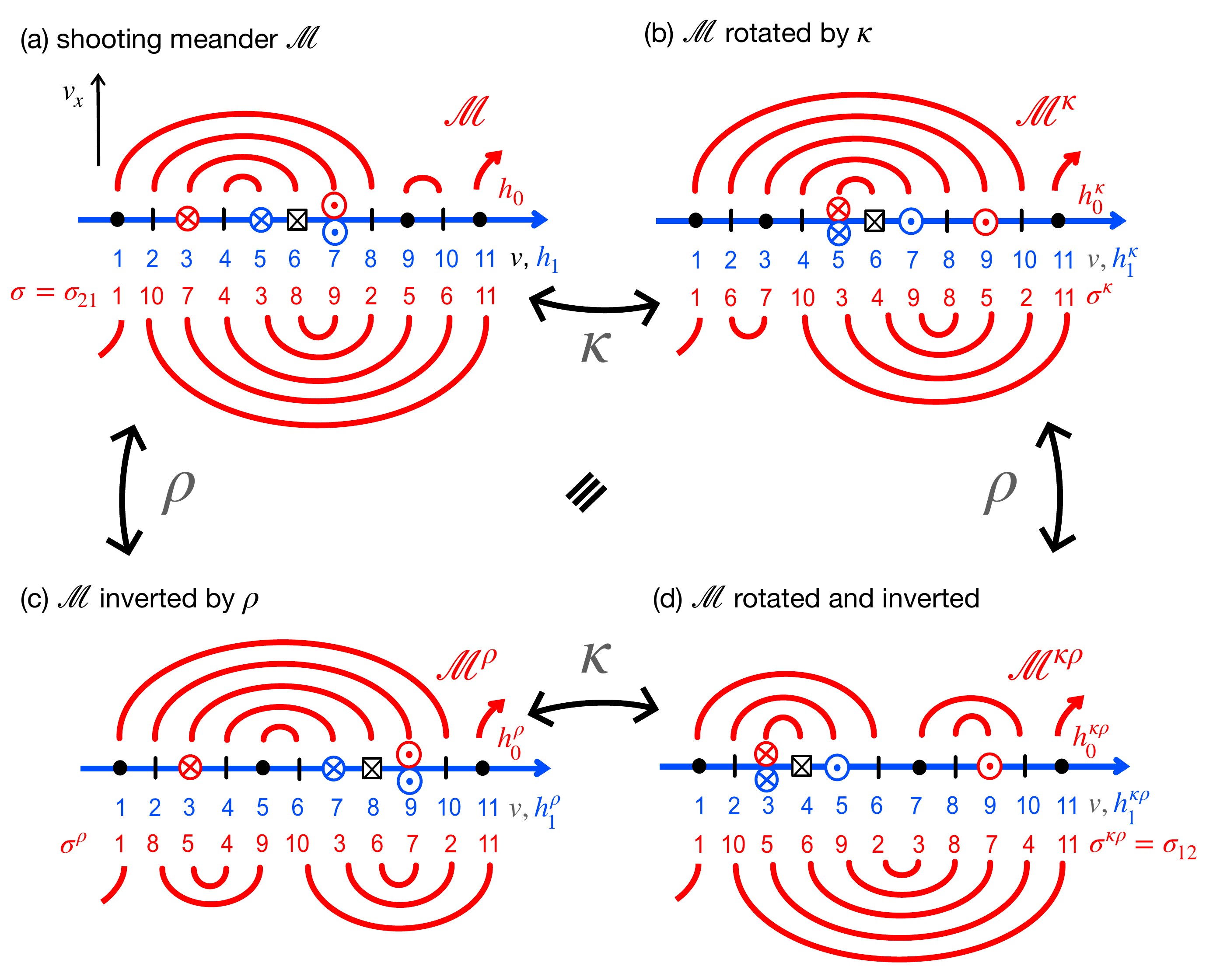}
\caption{\emph{
Example of four trivially equivalent 3-nose Sturm meanders $\mathcal{M}$, boundary orders $h_\iota$, and Sturm permutations $\sigma$.
A generalization of this diagram will be formulated explicitly, and proved, in lemma \ref{noseloclem}.
See also theorem \ref{krrhothm} and proposition \ref{sigmaprop}.
See fig.~\ref{3ball} for a legend of Morse indices.
(a) The Sturm meander $\mathcal{M}$ of the Sturm permutation $\sigma=(1\ 10\ 7\ 4\ 3\ 8\ 9\ 2\ 5\ 6\ 11)$. 
The three noses are located at horizontal positions $h_1=\{4,5\},\{6,7\}$, and $\{9,10\}$. 
Note $\mathcal{M}=\mathcal{M}_{41}$ and $\sigma=\sigma_{21}$, in the notation of definition \ref{Akrdefi}.
(b) The flip $\kappa$, applied to $\mathcal{M}$, produces the meander $\mathcal{M}^\kappa$, via rotation by $180^\circ$, the reverse boundary orders $h_\iota^\kappa=h_\iota \kappa$, and the conjugate permutation $\sigma^\kappa$.
See  \eqref{rotM} -- \eqref{rotsig}.
(c) The $x$-reversal $\rho$ replaces the Sturm permutation $\sigma$ by the inverse permutation $\sigma^\rho=\sigma^{-1}$, with the appropriate meander $\mathcal{M}^\rho$.
Note how the roles of $h_0$ and $h_1$ are interchanged.
For the equilibrium configuration $\mathcal{E}^\rho$, the Sturm attractor $\mathcal{A}^\rho$, and the connection graph $\mathcal{C}^\rho$ see fig.~\ref{3ball}.
(d) The combination $\kappa\rho$ of the two generators of trivial equivalences. 
Note the inversion $\rho$, from (a) to (c), and the commuting $180^\circ$ rotation $\kappa$, from (c) to (d).
}}
\label{equi}
\end{figure}

For example, and alternatively, let us describe the effect of $\kappa = -\text{id}$ on the meander $\mathcal{M}$, the boundary orders $h_\iota$\,, and on the Sturm permutations $\sigma$, algebraically.
The meander $\mathcal{M}$ is the (stylized) shooting image of the horizontal $v$-axis, in the $(v,v_x)$-plane, under the nonautonomous ODE flow \eqref{ODE}, evaluated from $x=0$ to $x=1$.
The involution \eqref{rot} therefore simply rotates $\mathcal{M}\subset \mathbb{R}^2$ by $180^\circ$, i.e. 
\begin{equation}
\label{rotM}
\mathcal{M}^\kappa:=-\mathcal{M}\,.
\end{equation}
The orientation of the meander curve, however, is reversed.
Abusing notation slightly, let $\kappa$ also denote the flip permutation
	\begin{equation}\label{flip}
	\kappa (j) := N+1-j
	\end{equation}
on $j \in \{ 1, \ldots, N\}$.
Then $\kappa$ reverses the boundary orders of the equilibria $\mathcal{E}^\kappa:= -\mathcal{E}$, at $x=\iota \in \{0,1\}$, respectively, i.e.
	\begin{equation}\label{roth}
	h_\iota^\kappa = \kappa h_\iota \kappa: \ \{1,\ldots,N\}\longrightarrow \mathcal{E}^\kappa\,.
	\end{equation}
Here $\kappa$ in $\kappa h_\iota \kappa$ refers to \eqref{rot}, on the left, and to \eqref{flip} on the right.
Therefore $\sigma = h_0^{-1}\circ h_1$ from \eqref{permdef}, alias the meander rotation \eqref{rotM}, leads to conjugation 
	\begin{equation}\label{rotsig}
	\sigma^\kappa = \kappa \sigma \kappa
	\end{equation}	
by the flip \eqref{flip}.
See the horizontal pair (a),(b) of fig.~\ref{equi} for the effect of the rotation $\kappa$ on the meander $\mathcal{M}$ of the Sturm permutation $\sigma=(1\ 10\ 7\ 4\ 3\ 8\ 9\ 2\ 5\ 6\ 11)$.
Similarly, the horizontal pair (c),(d) is $\kappa$-related.

\begin{figure}[t!]
\centering \includegraphics[width=\textwidth]{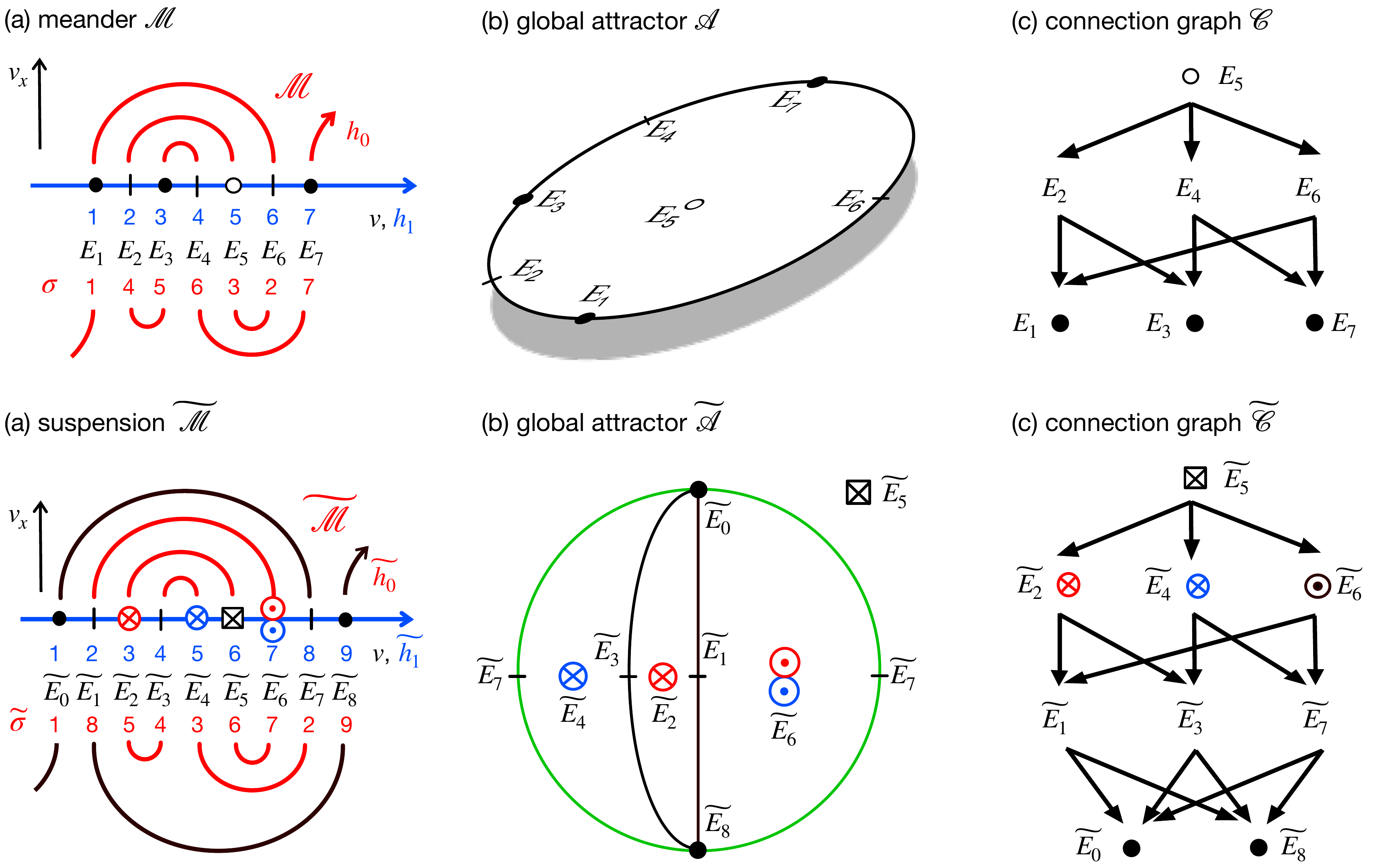}
\caption{\emph{
Suspensions $\widetilde{\mathcal{M}},\widetilde{\mathcal{A}},\widetilde{\mathcal{C}},\widetilde{\sigma}$, in the bottom row, of the 3-nose meander $\mathcal{M}$, the attractor $\mathcal{A}$, and the connection graph $\mathcal{C}$ for the Sturm permutation $\sigma=(1\ 4\ 5\ 6\ 3\ 2\ 7)$, top row.
See fig.~\ref{3ball} for a legend of Morse indices.
(a) The meander $\mathcal{M}$ of $\sigma$. 
The three noses are located at horizontal positions $h_1=\{2,3\},\{3,4\}$, and $\{5,6\}$. 
Equilibria $\mathcal{E}=\{E_1,\ldots,E_7\}$ are enumerated along the horizontal $h_1$-axis as $E_j=h_1(j)$.
(b) The 2-ball global attractor $\mathcal{A}= \mathrm{clos}\,W^u(E_5)$ is a triangle disk, topologically, with the stable equilibria $E_1,E_3,E_7$, as vertices, and the closures of the one-dimensional unstable manifolds of the saddles $E_2,E_4,E_6$, as sides.
(c) The connection graph $\mathcal{C}$ of the triangle attractor $\mathcal{A}$ from (b).
(d) The suspension meander $\widetilde{\mathcal{M}}$ of $\mathcal{M}$, with Sturm permutation $\widetilde{\sigma}=(1\ 8\ 5\ 4\ 3\ 6\ 7\ 2\ 9)$ and equilibria $\widetilde{E}_j:=\widetilde{h}_1(j+1),\ j=0,\ldots,8$.
See also proposition \ref{suspprop} and \eqref{lift},\eqref{suspapp}.
Black: added upper and lower arcs of $\widetilde{\mathcal{M}}$.
(e) The 3-ball global attractor $\widetilde{\mathcal{A}}$ of the suspension meander $\widetilde{\mathcal{M}}$ is the topological double cone suspension over the equatorially embedded repelling triangle disk $\mathcal{A}$. 
Indeed, the $1$-dimensionally unstable embedding $\mathcal{A} \subset \widetilde{\mathcal{A}}$ is based on the lifting  identification $E_j\mapsto\widetilde{E}_j$ of $\mathcal{E}\subset\widetilde{\mathcal{E}}$.
Note $i(\widetilde{E}_j)=i(E_j)+1$.
The cone vertices are the two added polar sinks $\widetilde{E}_0,\widetilde{E}_8$\, with Morse index $i=0$.
(f) The suspension connection graph $\widetilde{\mathcal{C}}$, correspondingly, contains a copy of $\mathcal{C}$, raised in Morse index $i$ by $1$, in the top three rows.
Towards each equilibrium $\widetilde{E}_0,\widetilde{E}_8$ in the polar $i=0$ bottom row, all $i=1$ saddles connect heteroclinically.
Transitivity in fact implies $\widetilde{E}_j \leadsto \widetilde{E}_0, \widetilde{E}_8$, for all lifted equilibria.
}}
\label{susp}
\end{figure}

Reversal $\rho$ of $x$, in contrast, interchanges the boundaries $x=\iota \in \{0,1\}$.
Therefore 
	\begin{equation}\label{invh}
	h_\iota^\rho = \rho h_{1-\iota}: \ \{1,\ldots,N\}\longrightarrow \rho\mathcal{E}\,,
	\end{equation}
and $\sigma = h_0^{-1}\circ h_1$ leads to inversion of the Sturm permutation
	\begin{equation}\label{invsig}
	\sigma^\rho = \sigma^{-1}\,.
	\end{equation}	
Graphically this amounts to pulling the original meander $h_0$ straight, and considering the original horizontal axis $h_1$ as a meander over the new horizontal axis (after an up-down reflection).
See the vertical pair (a),(c) of fig.~\ref{equi}, and also the pair (b),(d).
For a less trivial example see also proposition \ref{sigmaprop}.

A small subtlety arises, concerning isotropy $f^\gamma=f$ of nonlinearities under some trivial equivalence $\gamma \in \langle \kappa,\rho \rangle$.
Such $f$-isotropy implies permutation-isotropy $\sigma^\gamma = \sigma$, of course. 
However, we never proved the converse.
Although some nonlinearities $f$ will always realize isotropic permutations $\sigma=\sigma^\gamma$, by \cite{firo99}, we never proved realization by an $f$ with isotropy $\gamma$, i.e. such that $f=f^\gamma$.

To define the suspension of a dissipative meander $\mathcal{M}$, from $N$ to $N+2$ vertices, we first label the $N$ vertices $\mathcal{E}=\{E_1,\ldots,E_N\}$ of $\mathcal{M}$, along the horizontal axis, as $E_j:=h_1(j)$. 
See fig.~\ref{susp}(a) for a Sturm example with $N=7$ original equilibria.
For the $N$ vertices $\{\widetilde{E}_1,\ldots,\widetilde{E}_{N}\}$ among the suspension vertices $\widetilde{\mathcal{E}}=\{\widetilde{E}_0,\ldots,\widetilde{E}_{N+1}\}$, we choose a corresponding enumeration $\widetilde{E}_j=\widetilde{h}_1(j+1)$, for $j=1,\ldots,N$.
This embeds old vertices $\mathcal{E}\subset\widetilde{\mathcal{E}}$ into the suspension via the lifting identification
\begin{equation}\label{Elift}
E_j \mapsto \widetilde{E}_j\,.
\end{equation}
Henceforth, we write $E_j=\widetilde{E}_j$ under this identification.

We define the \emph{suspension} $\widetilde{\mathcal{M}}$ as an augmentation of $\mathcal{M}$ by two overarching arcs (black in fig.~\ref{susp}(d)): an upper arc from the first new vertex $\widetilde{E}_0$ to the last old vertex $\widetilde{E}_N=E_N$, and a lower arc from the first old vertex $\widetilde{E}_1=E_1$ to the last new vertex $\widetilde{E}_{N+1}$. 
This extends the previous definition of $\widetilde{h}_1$ to $\widetilde{h}_1(j):=\widetilde{E}_{j-1}$ for $j=1,\ldots,N+2$.

By construction, the number of meander-noses is invariant under suspension, for $N\geq 3$.
In the Sturm case, i.e. if our dissipative meanders $\widetilde{\mathcal{M}}$ are also Morse, our definition also extends to define the suspensions $\widetilde{\mathcal{A}}$ and $\widetilde{\mathcal{C}}$ of their attractors $\mathcal{A}$ and connection graphs $\mathcal{C}$.

More abstractly, however, our definition of suspension generalizes to dissipative meanders $\mathcal{M}$, which are not necessarily Sturm.
Indeed they may violate the Morse property $i_j\geq 0$ and hence may also violate $z_{jk}\geq 0$.
Abstractly however, dissipative meanders $\mathcal{M}$ still determine their permutations $\sigma$, Morse numbers $i_j$ and zero numbers $z_{jk}$ via \eqref{i},\eqref{z} -- even when those numbers lack any  ODE or PDE interpretation.
Sturm ``attractors'' $\mathcal{A}$ with actual ``equilibria'' and actual ``heteroclinic orbits'' $v_1\leadsto v_2$ cannot exist, of course, once negative ``Morse indices'' $i_j$ are involved.
By $z$-adjacency \eqref{wolfrum}, and blocking \eqref{block}, however, we can still define connection graphs $\mathcal{C}$.
Quite radically, indeed, we abuse the notation $i(E_j)\,, \ z(E_j-E_k)$ here, and even $E_j\leadsto E_k$\,, to denote the recursively defined quantities $i_j\,, \ z_{jk}$\,, and the relation $\leadsto$ defined abstractly via \eqref{wolfrum},\eqref{block}.
In particular our definition of meander suspensions readily extends to define the suspensions $\widetilde{\mathcal{M}}, \widetilde{\sigma}$, and $\widetilde{\mathcal{C}}$, even in non-Morse cases.
Of course, these remarks also extend the notions of trivial equivalences to merely dissipative non-Morse meanders, algebraically, by \eqref{flip},\eqref{rotsig},\eqref{invsig} instead of the explicit maps \eqref{rot},\eqref{inv}.

The following proposition justifies the name ``suspension''.
Indeed, we may view the suspension $\widetilde{\mathcal{A}}$ of a global Sturm attractor $\mathcal{A}$ as the double cone suspension of $\mathcal{A}$ itself, with respect to the two added polar cone vertices $\widetilde{E}_0$ and $\widetilde{E}_{N+1}$.
See fig.~\ref{susp} again.

\begin{prop}\label{suspprop}
For dissipative, but not necessarily Morse, meanders $\mathcal{M}$, the suspension defined above has the following properties, for all $1\leq j,k\leq N,\ j\neq k$: 
\begin{enumerate}[(i)]
  \item $\widetilde{\sigma}(1)=1$ and $\widetilde{\sigma}(N+2)=N+2$;
  \item $\widetilde{\sigma}(j+1) = N+2 - \sigma(j)=\kappa \sigma(j)+1$;
  \item $i(\widetilde{E}_0)=i(\widetilde{E}_{N+1})=0$;
  \item $i(\widetilde{E}_j)=i(E_j)+1$;
  \item $z(\widetilde{E}_j-\widetilde{E}_0) = z(\widetilde{E}_j-\widetilde{E}_{N+1}) = 0$;
  \item $z(\widetilde{E}_j-\widetilde{E}_k) = z(E_j-E_k)+1$;
  \item $\widetilde{E}_j \leadsto \widetilde{E}_k \quad \Longleftrightarrow \quad E_j \leadsto E_k$\,;
  \item $\widetilde{E}_j \leadsto \widetilde{E}_0\,, \widetilde{E}_{N+1}$\,, in case all $i_j\geq 0$.
\end{enumerate}
\end{prop}

\begin{proof}
Consider suspensions $\widetilde{h}_\iota: \{1,\ldots,N+2\}\rightarrow \widetilde{\mathcal{E}}$ and $\widetilde{\sigma}\in S_{N+2}$ of abstract ``boundary orders'' $h_\iota: \{1,\ldots,N\}\rightarrow\mathcal{E}$ which fix $1$ as well as $N$.
Define the dissipative meander permutation $\widetilde{\sigma}=\widetilde{h}_0^{-1} \widetilde{h}_1 \in S_N$, as in \eqref{permdef}.

Claim (i) then holds by construction.
To prove claim (ii), first note that $\widetilde{h}_1(j+1)=\widetilde{E}_j=E_j=h_1(j)$. 
Since the orders $\widetilde{h}_0$ and $h_0$ follow the shared part of the meanders $\widetilde{\mathcal{M}}$ and $\mathcal{M}$, in opposite directions, we also have $k+\widetilde{k}=N+2$ for $h_0(k):=E_j$ and $\widetilde{h}_0(\widetilde{k}):=\widetilde{E}_j$\,.
Together this proves (ii), if we substitute the flip $\kappa$ from \eqref{flip}.
Properties (iii)--(vi) can be derived from the explicit recursions \eqref{i} and \eqref{z}.
In particular, (iv) enters in (vi) via the term $z_{kk}$ which gets raised by 1 after suspension.

Property (vii) follows from Wolfrum blocking \eqref{wolfrum},\eqref{block}.
Indeed, (vi) implies that blocking \eqref{block} between lifted old vertices $v_1,v_2 \in \mathcal{E}$ by any new vertex $\widetilde{w} \in \{\widetilde{E}_0,\widetilde{E}_{N+1}\}$ cannot occur, because the $\widetilde{h}_\iota$-position of those new vertices is extremal and never between $v_1,v_2$.
By (vi), in contrast, any old blocking remains in effect.
This proves claim (vii).

In claim (viii) we assume $\mathcal{M}$ to be Morse, and hence Sturm.
In particular, this implies $z_{jk}\geq 0$, for all zero numbers.
Therefore (v),(vi) prevent blocking \eqref{block}, and (viii) follows from (iii),(iv) with \eqref{wolfrum}.
\end{proof}

\begin{cor}\label{suspcor}
For Sturm meanders $\mathcal{M}$ the following holds true.
\begin{enumerate}[(i)]
\item The suspension $\widetilde{\sigma}\in S_{N+2}$ of any Sturm permutation $\sigma\in S_N$ is Sturm.
\item All $i=1$ equilibria connect heteroclinically, in $\widetilde{\mathcal{C}}$, towards the two polar $i=0$ sinks $\widetilde{E}_0,\widetilde{E}_{N+1}$ in the bottom row.
\item The connection graph $\widetilde{\mathcal{C}}$ of the suspension contains the connection graph $\mathcal{C}$, lifted to the rows $i\geq 1$.
\end{enumerate}
\end{cor}

\begin{proof}
Claim (i) follows from proposition \ref{suspprop} (iii),(iv).
With (viii), this also proves claim (ii).
Claim (iii) then follows from (vii).
\end{proof}

We conclude with an elementary remark on the compatibility of our trivial equivalences $\kappa,\rho$ with suspensions. 
Proposition \ref{suspprop}(ii) easily implies that suspension preserves the notion of trivial equivalence, but not the particular equivalences by $\rho$ or $\kappa\rho$ :
\begin{align}
\label{suspk}
   \kappa: \qquad \widetilde{\kappa\sigma\kappa} &=  \kappa\widetilde{\sigma}\kappa\,;  \\
\label{suspr}
   \rho: \,\qquad \widetilde{\sigma^{-1}} &=  (\kappa\widetilde{\sigma}\kappa)^{-1}\,.
\end{align}
This works for all dissipative meander permutations, and is not restricted to the Morse case.
We have taken license here to denote the flip \eqref{flip} in, both, $S_N$ and $S_{N+2}$ by the same letter $\kappa$.

In the Sturm case, the realization of suspensions by nonlinearities $f(x,u,p)$ may be of applied interest in design.
For example, we may append a region $x\in[1,2]$ to the $x$-domain of \eqref{ODE},\eqref{PDE}.
Then suspension can be effected, in terms of $x$-profiles of equilibria like fig.~\ref{3ball}(b), if $f$ reverses the order of equilibria at the right boundary, as $x$ increases from $x=1$ to $x=2$.
This agrees well with proposition \ref{suspprop}(vi). 
Dissipativeness, of course, will require two new equilibria, e.g. homogeneous throughout $0\leq x\leq 2$: one at the top, and one at the bottom.

We formalize this construction as follows. 
Let $\Lambda: S_N \rightarrow S_{N+2}$ lift permutations by
\begin{equation}\label{lift}
(\Lambda\sigma)(j) \ :=\ 
\begin{cases}
      \sigma(j-1)+1,&\textrm{for } 2\leq j \leq N+1, \\
     j, &\textrm{for } j \in \{1,N+2\}.
\end{cases}
\end{equation}
Although $\Lambda$ does not preserve the Sturm property, corollary \ref{suspcor}(i) asserts that suspension does.
But we can now rewrite proposition \ref{suspprop}(i),(ii) as
\begin{equation}\label{suspapp}
\widetilde{\sigma} = \Lambda(\kappa\sigma)
\end{equation}
Indeed, this is the append-construction just described via the order reversing involution $\kappa\in S_N$ of \eqref{flip}.

Now suppose we prepend order reversion of the equilibria of $\sigma$ on an interval $x\in[-1,0]$, instead, and then lift by $\Lambda$ again.
We claim that this prepend-construction is another realization of the same suspension, 
\begin{equation}\label{suspprep}
\widetilde{\sigma} = \Lambda(\sigma\kappa),
\end{equation}
provided that flip isotropy $\sigma^\kappa=\kappa\sigma\kappa=\sigma$ holds for the original Sturm permutation $\sigma$, i.e. for $x\in[0,1]$.
Indeed, the lift \eqref{lift} commutes with inversion.
Therefore \eqref{suspapp},\eqref{suspr},\eqref{suspk} successively imply
\begin{equation}\label{prepapp}
\Lambda(\sigma\kappa) = \big(\Lambda((\sigma\kappa)^{-1})\big)^{-1} =\big(\Lambda(\kappa\sigma^{-1})\big)^{-1} = \big(\widetilde{\sigma^{-1}}\big)^{-1}=\kappa\widetilde{\sigma}\kappa = \widetilde{\kappa\sigma\kappa} = \widetilde{\sigma}.
\end{equation}
As a corollary, we conclude from \eqref{suspapp},\eqref{suspprep} that the spatial order of points $x$, where equilibrium profiles cross each other, may differ widely for one and the same Sturm permutation $\widetilde{\sigma}$.
It is instructive to compare the append- and prepend-constructions in the explicit cubic Chafee-Infante case $\sigma=\sigma_d$ of \eqref{cisigma} below; see also fig.~\ref{ci3}(b).

\section{Two noses: the Chafee-Infante paradigm}\label{ChIn}

In this section we study the sequence $\mathcal{M}_d\,,\ d\geq1$, of Sturm meanders with two noses and $2d$ arcs.
We first proceed completely abstractly, without reference to any specific nonlinearities, to derive the Sturm meanders $\mathcal{M}_d$\,, their Sturm global attractors $\mathcal{A}_d$\,, their connection graphs $\mathcal{C}_d$\,, and reversibility.
Instead of explicit calculations involving a specific nonlinearity, or simulations of mere anecdotal relevance, we exclusively rely on the general principles and concepts outlined in the previous sections.
Only as an after\-thought, we return to PDE \eqref{PDE} with symmetric cubic $f$ and parameter $\lambda$,
\begin{equation}\label{cubic}
f=\lambda^2 u (1-u^2), \qquad(d-1)\pi<\lambda<d\pi,
\end{equation}
as studied by Chafee and Infante \cite{chin74}.
Via the abstract 2-nose Sturm meander $\mathcal{M}_d$\,, we will see how our abstract global attractor $\mathcal{A}_d$ is actually orbit equivalent to the Chafee-Infante attractor of that explicit original example.

\begin{figure}[t]
\begin{center}
\centering \includegraphics[width=\textwidth]{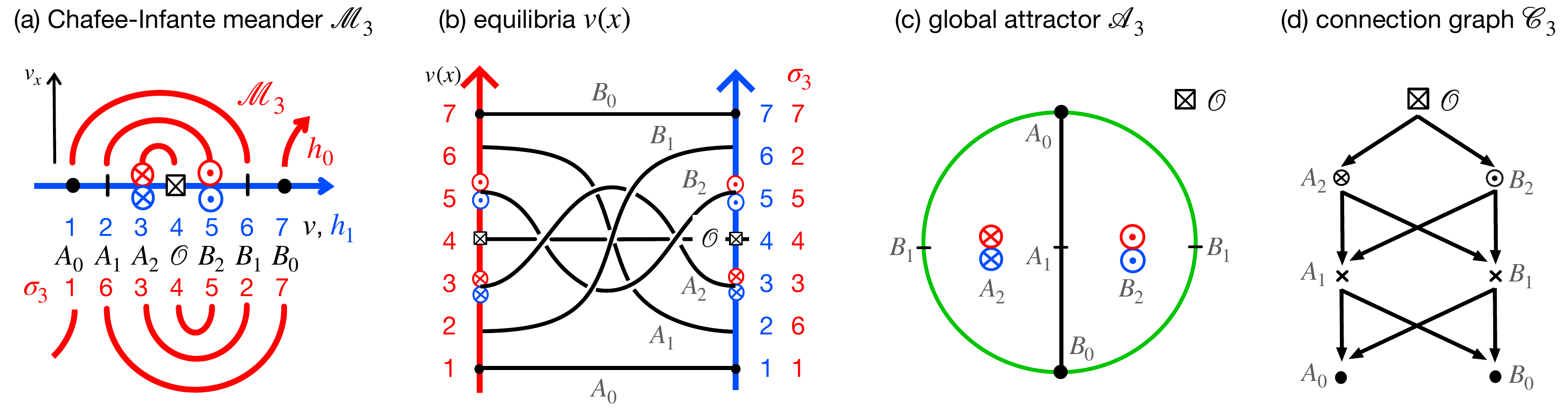}
\caption{\emph{
The Chafee-Infante attractor $\mathcal{A}_d$ of dimension $d=3$; 
see PDE \eqref{PDE} with cubic nonlinearity $f=\lambda^2 u(1-u^2)$, and parameter $2\pi<\lambda<3\pi$ according to \eqref{cubic}.
The four diagrams are analogous to fig.~\ref{3ball}(a)-(d).
(a) The associated schematic Sturm meander $\mathcal{M}_3$ has two noses, at equilibria $\{A_2,\mathcal{O}\}$ and at $\{\mathcal{O},B_2\}$. 
Note how the Sturm involution $\sigma_3$ is invariant under the trivial equivalence $\kappa$, because $f$ is odd, and under $\rho$, because $f=f(u)$.
(b) The seven equilibria $v(x)$ of the Chafee-Infante attractor $\mathcal{A}_3$\,.
The trivial equivalence $\kappa: v\mapsto -v$ interchanges each $A_j$ with $B_j$\,.
The trivial equivalence $\rho: x\mapsto 1-x$, in contrast, only swaps $A_1$ and $B_1$\,.
(c) The 3-ball Chafee-Infante attractor $\mathcal{A}_3=\mathrm{clos}\,W^u(\mathcal{O})$. 
Note time reversibility of the flow on its boundary 2-sphere $\partial W^u(\mathcal{O}) = \Sigma^2$.
(d) In the connection di-graph $\mathcal{C}_3$, ranked by Morse index, the trivial equivalences appear as automorphisms. Reversibility on the boundary 2-sphere $\Sigma^2$ below $\mathcal{O}$, in contrast, appears as a partial automorphism under the $180^\circ$ rotation $A_j \leftrightarrow B_{2-j}\,$.
It reverses the heteroclinic orientation of the di-graph $\mathcal{C}_3$, below $\mathcal{O}$, and the Morse ranking.
}}
\label{ci3}
\end{center}
\end{figure}

To pursue this program, let us start from just $d$ upper arcs, separately and without meanders in mind as yet.
Equivalently, the arcs define a balanced structure of $d$ pairs of opening and closing parentheses, ``$($'' and ``$)$'', also know as \emph{Dyck words} of length $2d$, as counted by the Catalan numbers.
For a historical reference see the habilitation thesis by Dyck on the word problem in combinatorial group theory \cite{dyck}.
Upper noses correspond to innermost pairs ``$()$''.
Any nonempty Dyck word has to contain at least one nose.
If the Dyck word only contains a single nose, then all parenthesis pairs, alias arcs, must be nested.
In section \ref{Back}, we already called such a total nesting a \emph{rainbow}.
Proceeding for lower arcs, analogously, we obtain another rainbow of $d$ nested lower arcs.
Dissipativeness requires the lower rainbow to be shifted one entry to the right, with respect to the upper rainbow. See fig.~\ref{ci3}(a).
Joining the two rainbows defines a unique double spiral which, automatically, turns out to be a dissipative meander $\mathcal{M}_d$\,, for any $d\geq 1$.
By construction, $\mathcal{M}_d$ possesses $d$ upper and $d$ lower arcs, each, over its $N=2d+1$ intersections with the horizontal axis.
Alas, we do not know yet whether $\mathcal{M}_d$ is Morse, and therefore Sturm.

Let us examine the associated meander permutation $\sigma_d$\,.
The numbers of arcs and noses are invariant under trivial equivalences \eqref{rotsig},\eqref{invsig}.
Therefore the unique dissipative 2-nose meander $\mathcal{M}_d$ of $2d$ arcs, and its meander permutation $\sigma_d$\,, are invariant under trivial equivalences:
\begin{equation}\label{ciequiv}
\sigma_d^\kappa = \sigma_d^\rho = \sigma_d\,.
\end{equation}
For $d=3$ see fig.~\ref{ci3}.
Suspension leaves the number of noses invariant, likewise, but increases the number $d$ of upper and of lower arcs by $1$, each. This proves
\begin{equation}\label{cisusp}
\widetilde{\,\sigma_d\,}=\sigma_{d+1}\,,
\end{equation}
notably without any calculation.
The case $\sigma_1 = (1\ 2\ 3)$ of $d=1$ is trivial.
In view of proposition \ref{suspprop}(iii),(iv), inductively, the dissipative meander $\mathcal{M}_d$ is therefore Morse, of maximal Morse index $d$, and hence Sturm.
We call the (stylized) Sturm meander $\mathcal{M}_d$ the \emph{Chafee-Infante meander of dimension} $d$.

To determine the connection graph $\mathcal{C}_d$, we relabel the $N=2d+1$ equilibria $E_j$ along the horizontal axis:
\begin{equation}\label{ciequi}
E_j\ =\ h_1(j) \ =:\ 
\begin{cases}
      \ A_{j-1}\,, & \textrm{for } 1\leq j < d+1, \\
      \ \mathcal{O}, & \textrm{for } j=d+1, \\
      \ B_{N-j}\,, & \textrm{for } d+1 < j \leq N.
\end{cases}
\end{equation}
Starting from the trivial $1$-dimensional Sturm case $d=1$ with $\mathcal{O}\leadsto \{A_0,B_0\}$, we arrive at the general case by successive suspension \eqref{cisusp}. 
In view of suspension corollary \ref{suspcor}, this identifies the connection graph $\mathcal{C}_d\,,\ d\geq1,$ to be given by
\begin{align}
\label{CdO}
    \mathcal{O}&\leadsto \{A_{d-1},B_{d-1}\}, \textrm{ and}   \\
\label{CdAB}
    A_j,\,B_j&\leadsto \{A_{j-1},B_{j-1}\}, \textrm{ for } 1\leq j < d. 
\end{align}
See fig.~\ref{ci3}(d).
In particular, suspension \eqref{cisusp} or just cascading identify the Morse indices
\begin{equation}\label{cimorse}
i(A_j)=i(B_j)=j, \qquad i(\mathcal{O})=d=\dim \mathcal{A}_d\,.
\end{equation}
By transitivity, $\mathcal{O}$ connects to any other equilibrium, heteroclinically.
Therefore, $\mathcal{A}_d= \mathrm{clos}\, W^u(\mathcal{O})$ is a Sturm $d$-ball and $\partial W^u(\mathcal{O}) = \Sigma^{d-1}$ is a ($d$-1)-sphere.
By \eqref{CdAB}, the connection graph $\mathcal{C}_d$ is reversible on $\Sigma^{d-1}$, e.g.~under the involutive reversor
\begin{equation}\label{cirev}
\mathcal{R}A_j \longleftrightarrow B_{d-1-j}\,,
\end{equation}
for $0\leq j<d$; see \eqref{rev}.

Although we did not use this above, at all, we at least mention that the Sturm permutation $\sigma_d$ of $\mathcal{M}_d$\,, i.e.~with $N=2d+1$ intersections, is given explicitly by 
\begin{equation}\label{cisigma}
\sigma_d(j)  \ =\  
\begin{cases}
      \  j, &\textrm{ for odd } j,  \\
      \ N+1-j, &\textrm{ for even } j.
\end{cases} 
\end{equation}
Indeed, the nested arcs are obtained correctly by the constant sum of the horizontal positions $\sigma^{-1}(j)$ of their successive endpoints, along the meander:
\begin{equation}\label{cinest}
\scalebox{0.87}[0.95]
	{$\sigma^{-1}(j)+\sigma^{-1}(j+1)=\sigma(j)+\sigma(j+1)=
\begin{cases}
      \  j+(N+1-(j+1))=N, &\textrm{ for odd } j,  \\
      \ (N+1-j) + (j+1)=N+2, &\textrm{ for even } j.
\end{cases}
$}
\end{equation}
Of course, our claims \eqref{ciequiv}--\eqref{cirev} could also have been derived from the explicit form \eqref{cisigma}, directly via \eqref{i}--\eqref{block} -- and without any deeper understanding.

We have already mentioned that the 1974 Chafee-Infante version $\mathcal{A}_d$ of \eqref{PDE} had been studied for the cubic nonlinearity $f(u)=\lambda^2 u(1-u^2)$, originally, albeit under Dirichlet boundary conditions; see \eqref{cubic} and the original paper \cite{chin74}.
Their method was local bifurcation analysis of the trivial equilibrium $\mathcal{O}: v\equiv 0$.
Note $i(\mathcal{O})=d$, for $(d-1)\pi<\lambda<d\pi$, under Neumann boundary conditions, by elementary linearization.
The second order ODE \eqref{ODE} is Hamiltonian integrable, for nonlinearities $f=f(u)$.
For the hard spring cubic nonlinearities $f$, the minimal periods $p(a)$ of $v(x)$ at $\lambda=1$ grow monotonically with their amplitude $a=v(0)\in (0,1)$ at $v_x(0)=0$. 
Note the limit $p(0)=2\pi$.
Rescaling $x$ as in $v(\lambda x)$, we see that $v$ reappears as a rescaled solution at $\lambda=n p(a)/2$, for any nonzero integer $n$.
See fig.~\ref{ci3}(b) for such rescaled equilibrium profiles in case $d=3$.
In particular, this produces a (stylized) shooting meander which, by monotonicity of the periods, coincides with the  Sturm meander $\mathcal{M}_d$\,, and hence determines the Sturm permutation $\sigma_d$ of \eqref{cisigma}.

For an early geometric description of the Chafee-Infante attractor $\mathcal{A}_d$\,, for low dimensions  $d=1,2,3$, see section 5.3 of \cite{he81}.
In 1985, Henry achieved the first description of $\mathcal{A}_d$ for general $d$ \cite{he85}.
His description was based on a nodal property akin to \eqref{zdrop}, and on a careful geometric analysis of unstable and center manifolds at the sequence of pitchfork bifurcations from $\mathcal{O}: v\equiv 0$, at $\lambda=n\pi$.
See fig.~\ref{ci3}(c).

In section 5 of \cite{firo20}, we have discussed the Sturm complex of the Chafee-Infante attractors in the more refined setting of signed hemisphere decompositions, which also leads to fig.~\ref{ci3}(c).
This also provides extremal characterizations of the Chafee-Infante attractor $\mathcal{A}_d$, among all Sturm attractors:
\begin{description}
  \item[max:] \emph{Among all Sturm attractors with $N=2d+1\geq 3$ equilibria, $\mathcal{A}_d$ is the unique Sturm attractor with the maximal possible dimension, $d$.}
  \item[min:] \emph{Among all Sturm attractors of dimension $d\geq 1$, $\mathcal{A}_d$ is the unique Sturm attractor with the smallest possible number of equilibria, $N=2d+1$.}
\end{description}
The two claims follow, e.g., from the connection graph.
In fact, each unstable hyperbolic equilibrium $v$ must connect, heteroclinically, to at least two other equilibria $v_\pm$, such that $z(v_\pm-v)=i(v)-1$ and $v_-<v<v_+$ at $x=0$. See also \cite{fi94}.

Topological Conley index and the connection matrix have been employed by Mischaikow \cite{mi95}, to establish heteroclinic orbits in larger classes of gradient-like PDEs with equilibrium configurations of Chafee-Infante type.
This technique establishes the existence of some (possibly non-unique) heteroclinic orbit $u(t,x)$ between the sets $\{A_j,B_j\}$ and $\{A_{j-1},B_{j-1}\}$\,.
Acting on $u(t,x)$ with the Klein 4-group of symmetries $\langle \kappa, \rho\rangle$, generated by \eqref{rot},\eqref{inv} in the Sturm setting, we obtain the four required heteroclinic orbits \eqref{CdAB}.
Indeed $\kappa$, alias $-\mathrm{id}$, interchanges each $A_j$ with $B_j$; see \eqref{rot} and fig.~\ref{ci3}(b). Inversion $\rho$, in contrast, performs the same interchange for odd $j$, only; see \eqref{inv}. 
Since the Morse levels $j$ and $j-1$ are of opposite even/odd parity, this generates the four required heteroclinic orbits.
The argument for the heteroclinic orbits emanating from the equilibrium $\mathcal{O}$, of top Morse index, is analogous.

From an applied point of view, \cite{mi95} greatly extends the Chafee-Infante paradigm beyond the requirement of Sturm zero numbers -- as long as a variational structure remains intact, with the same (minimal) configuration of equilibria, symmetries, and Morse indices.
This includes damped wave equations and other applications. See also \cite{hami91}.

An explicit ODE template for Chafee-Infante attractors has also been discussed in \cite{mi95}.
Consider the $d\times d$ diagonal matrix $\mathbf{Q} e_j = \mu_j e_j\,,\ j=0,\ldots, d-1$, with $d$ simple real eigenvalues $\mu_0>\ldots>\mu_{d-1}$\,.
Let $\mathbf{P}_\Phi:= (\mathrm{id} - \Phi \Phi^T)$ denote orthogonal projection onto the tangent space of the unit sphere at $\Phi\in S^{d-1} \subset \mathbb{R}^d$.
In polar coordinates $(\mathbf{r},\Phi)$, define the flow on $\mathbb{R}^d$ by
\begin{equation}\label{cimi}
\begin{aligned}
    \dot{\mathbf{r}} &=\mathbf{r}(1-\mathbf{r}^2)\,;   \\
    \dot{\Phi} &=\mathbf{P}_\Phi \, \mathbf{Q}\Phi  \,.
\end{aligned}
\end{equation}
The global attractor of this flow is orbit equivalent to the Chafee-Infante flow on $\mathcal{A}_d$\,.
In fact, the antipodal pairs $B_j=-A_j$ of Chafee-Infante equilibria correspond to the unit eigenvectors $\pm e_j$ of $\mu_j$\,. 
The equilibrium $\mathcal{O}$ maps to the origin $\mathbf{r}=0$.
Note how the connection graph \eqref{CdO},\eqref{CdAB} of the Chafee-Infante flow is realized by \eqref{cimi}.

Later work in the Sturm context addressed general autonomous nonlinearities $f(u)$; see for example \cite{brfi88, brfi89, firowo11}.
The paradigm of pitchfork bifurcations has been beautifully extended by Karnauhova, with many pictures, in \cite{ka17}.
With the pitchforkable class essentially well-understood since \cite{he85}, however, the simplest non-pitchforkable example had been discovered in \cite{ro91}.
Since none of our 3-nose meanders of dimension three or higher will fall into the pitchforkable class, either, we have to take another approach instead.
We will progress further along the more promising abstract path which, as a warm-up, we have just sketched for the Chafee-Infante problem.

\section{Three noses: main results}\label{Res}
In this section we present our main results on meanders with three noses.
The Chafee-Infante case of only two noses, discussed in the previous section \ref{ChIn}, will serve as a paradigm not to be skipped.
The general case of Sturm meanders $\mathcal{M}_{pq}$ reduces to the sequences $p=r(q+1),\ r,q\geq1$; see theorems \ref{nonmorse}, \ref{krmuthm}.
As usual, these come with their entourage of Sturm permutations $\sigma_{rq}$\,, associated Sturm attractors $\mathcal{A}_{rq}$\,, and connection graphs $\mathcal{C}_{rq}$ (definition \ref{Akrdefi}).
In theorem \ref{krmuthm}, we determine the Morse polynomial, i.e. we count the number of equilibria for each Morse index.
The Morse polynomials of $\mathcal{A}_{rq}$ and $\mathcal{A}_{qr}$ coincide; see corollary \ref{krmucor}.
In fact, the Sturm attractors $\mathcal{A}_{rq}$ turn out to be trivially equivalent to $\mathcal{A}_{qr}$\,, by theorem \ref{krrhothm} and corollary \ref{krrhocor}.
Geometrically, these are Sturm balls $\mathcal{A}_{rq}=\mathrm{clos}\,W^u(\mathcal{O})$ of dimension $i(\mathcal{O})=r+q$ (theorem \ref{krballthm}).
Finally, theorem \ref{krrevthm} asserts that the connection graph remains time reversible on the invariant boundary sphere $\partial \mathcal{A}_{rq} = \Sigma^{r+q-1}$.

To not clutter our conceptual approach by baroque notation, we will refrain from proving our results in full generality.
Instead, we only address the simplest interesting case $r=1$, i.e.~$\mathcal{A}_{1q}$ and $\mathcal{A}_{q1}$\,, in section \ref{r=1}.
For $r>1$ see our sequel \cite{firo23}.

As in the Chafee-Infante case of just two noses, section \ref{ChIn}, we start from upper and lower Dyck words, i.e. from pairs ``$( \ldots)$'' of opening and closing parentheses.
Three noses, this time, correspond to three innermost pairs ``$()$'' with their associated nestings.
Rotating by trivial equivalence $\kappa$, we may assume two nests to be upper, and one lower.
In other words, lower arcs form a rainbow, as before.
The upper Dyck word, however, takes the general form 
\begin{equation}\label{upperD}
(^s\,(^p)^p\,(^q)^q\,)^s,
\end{equation}
where exponents indicate repeated parentheses.
In the same notation, the lower rainbow becomes $(^{s+p+q})^{s+p+q}$, shifted by one vertex to the right.
Up to $s$ suspensions ``$\sim$'', and possibly a rotation $\kappa$ as in \eqref{flip}, we will therefore assume $s=0$.
Let $\mathcal{M}_{pq}$ denote the resulting dissipative configuration of upper arcs \eqref{upperD} and the lower $p+q$-rainbow.
Note how the special cases of either $p=0$ or $q=0$ are Chafee-Infante attractors, of two noses.

\begin{thm}\label{nonmorse}
With the above notation the following holds true for $r,p,q\geq1$.
\begin{enumerate}[(i)]
  \item $\mathcal{M}_{pq}$ is a dissipative meander if, and only if, $(p-1,q+1)$ are co-prime and $p\geq2$.
  \item For $p\neq r(q+1)$, any dissipative meander $\mathcal{M}_{pq}$ fails to be Morse.
\end{enumerate}
\end{thm}

\begin{figure}[t]
\begin{center}
\centering \includegraphics[width=\textwidth]{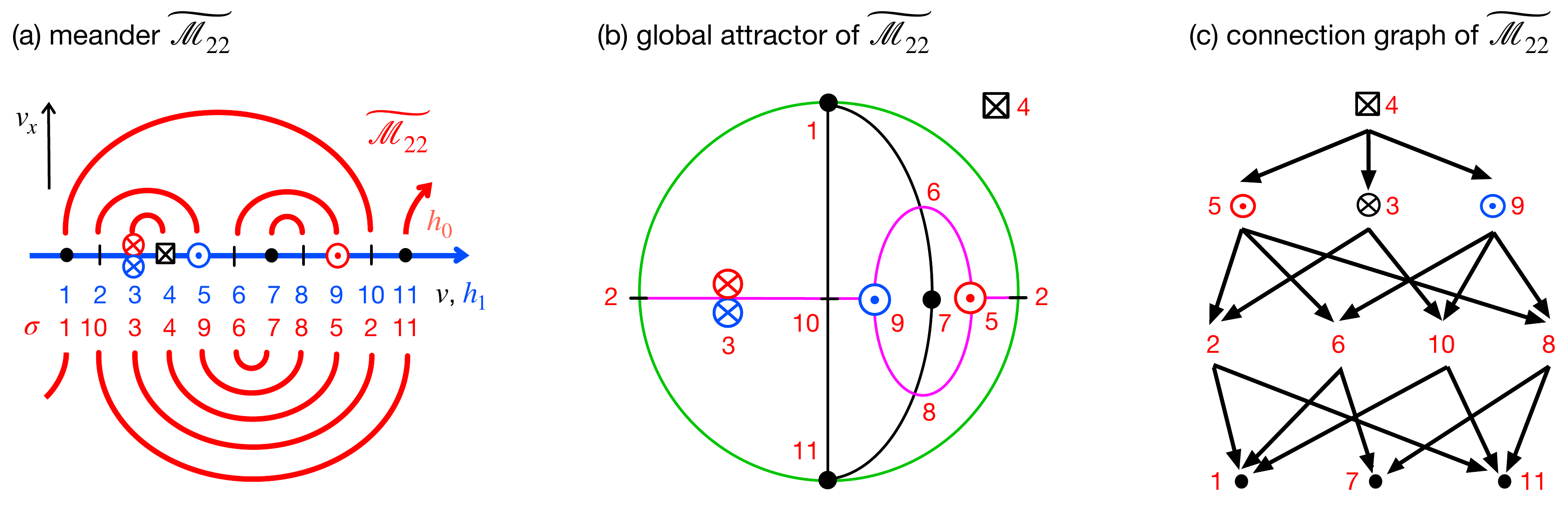}
\caption{\emph{
The suspension $\widetilde{\mathcal{M}_{22}}$ of the non-Morse, 3-nose meander $\mathcal{M}_{pq}$\,, with $p=q=2$.
The three diagrams are analogous to fig.~\ref{susp}(d)-(f).
Equilibria are enumerated as $\{1,\ldots,11\}$,  in red, along the meander path $h_0$\,.
(a) 
The suspended Sturm permutation $\sigma = (1\ 10\ 3\ 4\ 9\ 6\ 7\ 8\ 5\ 2\ 11)$ arises from the non-Morse meander permutation $(1\ 8\ 7\ 2\ 5\ 4\ 3\ 6\ 9)$.
(b) The 3-ball global attractor of $\widetilde{\mathcal{M}_{22}}$\,. Note time reversibility of the flow on the boundary 2-sphere.
(c) The connection di-graph of $\widetilde{\mathcal{M}_{22}}$\,, ranked by Morse index.
Time reversal on the boundary 2-sphere is induced by the reversor $automorphism \mathcal{R}$ which swaps the equilibria $(7,8,9,10,11) \leftrightarrow (3,2,1,6,5)$, in this order. This amounts to a rotation by $180^\circ$ of the lower three rows of the connection graph, and reversal of all arrows.
}}
\label{m22susp}
\end{center}
\end{figure}

We will prove theorem \ref{nonmorse} in section \ref{Nonmorse}.
Note that the non-Morse 3-nose cases (i), with $p\neq r(q+1)$ are not a lost cause, from the Sturm PDE point of view \eqref{PDE}.
Indeed, suspension proposition \ref{suspprop}(iv) always provides a minimal number $s=-\min i(E_j)\geq 0$ of suspensions after which $\mathcal{M}_{pq}$ becomes Morse, and hence Sturm.
See fig.~\ref{m22susp} for the non-Morse 3-nose example $p=q=2\neq r(q+1)$.
We will pursue those cases further in our sequel \cite{firo23}.

Let us now focus on the 3-nose cases $p=r(q+1)$, which are complementary to theorem \ref{nonmorse}(ii).
Then $p-1,q+1$ are automatically co-prime, because $r(q+1)-(p-1)=1$.
The following theorem shows that all cases $p=r(q+1)$ do lead to Morse meanders $\mathcal{M}_{pq}$\,, and therefore to Sturm attractors.
The rotation $\sigma:=\kappa\sigma_{11}\kappa$ of the simplest case $r=q=1$ has already served in fig.~\ref{susp}(a)-(c), to illustrate suspension.
We therefore assume $r,q\geq 1, \ rq>1$, for the rest of this paper.
Proofs of the next four theorems will be given in section \ref{r=1}, for the simplest interesting case $r=1$, only. See \cite{firo23} for general $r\geq 1$.

\begin{thm}\label{krmuthm}
For $p=r(q+1)$, let $m_{rq}(i)$ count the vertices with Morse number $i$, in the dissipative meander $\mathcal{M}_{pq}$\,. 
Then, for $r,q\geq 1, \ rq>1$, the nonzero Morse counts are given by
\begin{equation}\label{krmu}
m_{rq}(i)\ =\ 
\begin{cases}
      3+2i,& \textrm{for}\quad 0\leq i< \min\{r,q\};\\
      2+2\min\{r,q\},& \textrm{for}\quad \min\{r,q\}\leq i < \max\{r,q\};\\
      2(r+q)+1-2i,& \textrm{for}\quad \max\{r,q\}\leq i\leq r+q.
\end{cases}
\end{equation}
In particular, all such meanders $\mathcal{M}_{pq}$ are Sturm.
\end{thm}

\begin{cor}\label{krmucor}
The Morse count functions $i\mapsto m_{rq}(i)$ have the following symmetry properties.
\begin{enumerate}[(i)]
  \item  Up to ordering, the subscript set $\{r,q\}$ is determined by $m_{rq}$\,.
  \item Conversely, the subscript set determines $m_{rq}=m_{qr}$\,.
  \item For all $0\leq i < r+q$, we have $m_{rq}(i) = m_{rq}(r+q-1-i)$.
\end{enumerate}
\end{cor}

\begin{proof}
To prove (i), just note $\{r,q\}=\{\min\{r,q\},\max\{r,q\}\}$\,.
Claim (ii) follows from $r+q=\min\{r,q\}+\max\{r,q\}$\,.
To prove (iii), insert \eqref{krmu}.
\end{proof}

\begin{defi}\label{Akrdefi}
For any $r,q\geq1$, we call $\mathcal{M}_{r(q+1),q}$ a \emph{primitive} $3$\emph{-nose meander.}
For the Sturm entourage of $\mathcal{M}_{r(q+1),q}$\,, we denote the associated \emph{primitive}
Sturm permutation as $\sigma_{rq}$\,, the \emph{primitive} Sturm attractor as $\mathcal{A}_{rq}$\,, and the \emph{primitive} connection graph as $\mathcal{C}_{rq}$\,.
\end{defi}

\begin{thm}\label{krrhothm}
The primitive Sturm permutations $\sigma_{rq}$ and $\sigma_{qr}$ are trivially equivalent under the involutive product $\kappa\rho$ of \eqref{rotsig} and \eqref{invsig}, i.e.
\begin{equation}\label{krrho}
\sigma_{qr}=\kappa(\sigma_{rq})^{-1}\kappa.
\end{equation} 
\end{thm}

\begin{cor}\label{krrhocor}
The primitive 3-nose Sturm attractors $\mathcal{A}_{rq}$ and $\mathcal{A}_{r'q'}$ are orbit equivalent if, and only if, their subscript sets coincide, up to ordering.
In fact $\mathcal{A}_{rq}$ and $\mathcal{A}_{qr}$ are trivially equivalent, under the involutive product $\kappa\rho$ of \eqref{rot} and \eqref{inv}.
\end{cor}

\begin{proof}
Suppose $\mathcal{A}_{rq}$ and $\mathcal{A}_{r'q'}$ are orbit equivalent. Then their Morse counts coincide, and the first claim follows from corollary \ref{krmucor}.\newline
Conversely, suppose their subscript sets coincide, but with reversed order. Then the trivial equivalence of the attractors follows from theorem \ref{krrhothm} and section \ref{Rot}.
\end{proof}

The corollary contains an elementary hint why the dynamics on our 3-nose attractors have never been addressed in the literature, so far.
Indeed, trivial rotation equivalence \eqref{rot} switches rainbows of the 3-nose meanders between the lower and the upper side.
For $r\neq q$, this provides a total of four different meanders, under the Klein 4-group $\langle \kappa,\rho\rangle$ of trivial equivalences.
In particular, none of the trivial equivalences $\kappa,\rho,\kappa\rho$ can act as an isotropy on the nonlinearity $f$.
Therefore, the four related nonlinearities
\begin{equation}
\label{fiso}
f(x,u,p),\quad -f(x,-u,-p),\quad f(1-x,u,-p),\quad -f(1-x,-u,p)
\end{equation}
all have to be distinct functions, on any primitive (or suspended) global attractor $\mathcal{A}_{rq}$\,.
Outright, this excludes ODE-integrable nonlinearities $f=f(u)$, or $x$-reversible nonlinearities $f=f((x-\tfrac{1}{2})^2,u,p^2)$, as models.
Of course, the realization of Sturm meanders by ``certain'' dissipative nonlinearities $f$ is guaranteed by \cite{firo99}.
But the remaining non-integrable choices are so cumbersome to analyze, in any detail, that they have deterred all explicit efforts, so far.

In the ``symmetric'' case $r=q$, theorem \ref{krrhothm} reveals the only nontrivial isotropy $\kappa\rho$, in the Klein 4-group of trivial linear equivalences. 
In particular, the rainbow argument above shows that $f$ still cannot be $\kappa$-isotropic.
Admittedly, \eqref{krrho} suggests to study $f$ which commute with $\kappa\rho$, i.e. $f(x,u,p)=-f(1-x,-u,p)$.
However, $f=f(u)$ is still excluded, because $f(x,u,p)$ and $f(1-x,u,-p)$ must remain distinct.

Note the Morse count $m_{rq}(r+q)=1$ at maximal $i=r+q$; see \eqref{krmu}. 
Let $\mathcal{O}$ denote that unique equilibrium in $\mathcal{A}_{rq}$ of maximal Morse index $i(\mathcal{O})=r+q=\dim \mathcal{A}_{rq}$\,.

\begin{thm}\label{krballthm}
The primitive Sturm attractor $\mathcal{A}_{rq}$ is the closure of the unstable manifold of the single equilibrium $\mathcal{O}\in \mathcal{A}_{rq}$\,.
I.e., $\mathcal{A}_{rq}$ is a Sturm ball of dimension $r+q$.
\end{thm}

For even dimension $d=r+q$, corollary \ref{krmucor}(iii) makes it trivial to check that the Euler characteristic $\chi$ of the global attractor $\mathcal{A}_{rq}$ satisfies
\begin{equation}\label{euler}
\chi(\mathcal{A}_{rq}) := \sum_{i=0}^{d} (-1)^i m_{rq}(i) = m_{rq}(r+q) = 1,
\end{equation}
as is proper for any global attractor \cite{ha88}.
For odd $d$, this useful test of \eqref{krmu} is less trivial to check.
Taken $\mathrm{mod}\ 2$, of course, it again implies that the total number $N$ of equilibria must be odd.

\begin{cor}\label{krballcor}
With dimension replaced by $r+q+s,\ s>0$, theorem \ref{krballthm} remains valid for any $s$-fold suspension of $\mathcal{A}_{rq}$\,.
\end{cor}

\begin{proof}
By corollary \ref{suspcor}, suspensions of Sturm balls are Sturm balls.
\end{proof}

More surprisingly than in the Chafee-Infante case, we still observe time reversibility on the sphere boundary of the primitive 3-nose Sturm global attractors $\mathcal{A}_{rq}$ -- in spite of the parabolic, diffusion-dominated nature of the underlying original PDE \eqref{PDE}.

\begin{thm}\label{krrevthm}
The connection graph $\mathcal{C}_{rq}$ is reversible on the flow-invariant boundary sphere $\Sigma^{r+q-1}=\partial\mathcal{A}_{rq}=\partial W^u(\mathcal{O})$ of the primitive Sturm ball $\mathcal{A}_{rq} = \mathrm{clos}\ W^u(\mathcal{O})$.
\end{thm}

The reversibility on the boundary sphere $\Sigma^{r+q-1}$, of course, is a much deeper reason for the symmetry of the Morse count function $i\mapsto m_{rq}(i)$, for $0\leq i < r+q$, which we have already noticed in corollary \ref{krmucor}(iii).
Indeed, the reversor $\mathcal{R}$ on $\mathcal{C}_{rq}\setminus\mathcal{O}$ swaps equilibria of Morse indices $i$ and $r+q-1-i$.

The reversibility statement of theorem \ref{krrevthm} is violated for the $s$-fold suspension $\mathcal{A}$ of any primitive 3-nose attractor $\mathcal{A}_{rq}$\,.
Indeed, let $m_{rq}$ denote the Morse count \eqref{krmu} of $\mathcal{A}_{rq}$\,.
Then proposition \ref{suspprop}(iii),(iv) raises the Morse count $m$ of $\mathcal{A}$ to be
\begin{equation}\label{krmususp}
m(i)\ =\ 
\begin{cases}
      2,& \textrm{for}\quad 0\leq i< s,\\
      m_{rq}(i-s)& \textrm{for}\quad s\leq i \leq r+q+s,
\end{cases}
\end{equation}
for some $s>0$.
Therefore $\Sigma^{r+q+s-1}=\partial\mathcal{A}=\partial W^u(\mathcal{O})$ contains $m(0)=2$ sink equilibria, at Morse level $i=0$. 
At the highest Morse level $i=r+q+s-1$ in  $\Sigma^{r+q+s-1}$, in contrast, we encounter $m(r+q+s-1)=3$ equilibria.
This asymmetry violates reversibility.

\section{Non-Morse meanders with three noses}\label{Nonmorse}

In this section we prove theorem \ref{nonmorse}.

Claim (i) states that the dissipative arc configuration of $p$ nested upper arcs followed by $q$ nested upper arcs, and a right shifted lower $(p+q)$-rainbow, is a meander if, and only if, $p\geq 2$ and $(p-1,q+1)$ are co-prime.

The case $p=1$ is trivially discarded: all upper arcs of the nonempty $q$-nest close up to become circles, with the corresponding inner arcs of the lower rainbow. This contradicts the meander property.

For $p\geq2$, let us remove the outermost arc of the upper $p$-nest and, instead, stack it onto the upper $q$-nest. 
The resulting closed arc configuration now features upper nests of $p-1$ and $q+1$ arcs over the same lower rainbow.
This closing construction has been described and studied in \cite{fica13}, in terms of certain Cartesian billiards.
See also \cite{ka17,del18,zog20}, and the many references there.
The closing provides a closed Jordan curve if, and only if, the original dissipative arc configuration is a meander.
In other words, we obtain \emph{closed meanders} from dissipative meanders, and vice versa.

Let us now return to the dissipative arc configuration of $\mathcal{M}_{pq}$ with $p\geq 2$.
By (6.1) of \cite{fica13}, the greatest common divisor of $p-1$ and $q+1$ counts the connected components of the resulting closed arc configuration.
The proof was recursive, via the Euclidean algorithm for $(p-1,q+1)$.
This proves claim (i).

It remains to show, (ii), that the dissipative meander $\mathcal{M}_{pq}$ fails to be Morse, if $2\leq p\neq r(q+1)$ for any integer $r\geq 1$.


We first consider the case $2\leq p<q+1$.
We label equilibria such that $h_1=\mathrm{id}$.
Then $A:=2p+1$ and $B:=2p+2q$ are the left and right endpoints of the uppermost arc in the upper $q$-nest.
By \eqref{i}, Morse numbers of $h_\iota$-adjacent vertices are adjacent.
Obviously $B$ is $h_1$ adjacent to $N=2p+2q+1$.
By dissipativeness, $i(N)=0$.
Adjacency implies $i(B)=\pm 1$.
In case $i(B)=-1$, we are done.

In case $i(B)=+1$, we obtain $i(A)=0$ because the meander arc $AB$ turns left from $B$ to $A$; see \eqref{i} again.
Now consider the preceding lower rainbow arc from $B'$ to $A':=A=2p+1$.
Since $a+b=2p+2q+3$ for the two endpoints $a,b$ of any lower rainbow arc, our assumption $p<q+1$ implies $B'>A'$: the lower arc $A'B'$ turns left, from $A'$ to $B'$.
But we already know $i(A')=i(A)=0$.
Therefore \eqref{i} implies a negative Morse index $i(B')=-1$, and we are done again.

\begin{figure}[t]
\begin{center}
\centering \includegraphics[width=\textwidth]{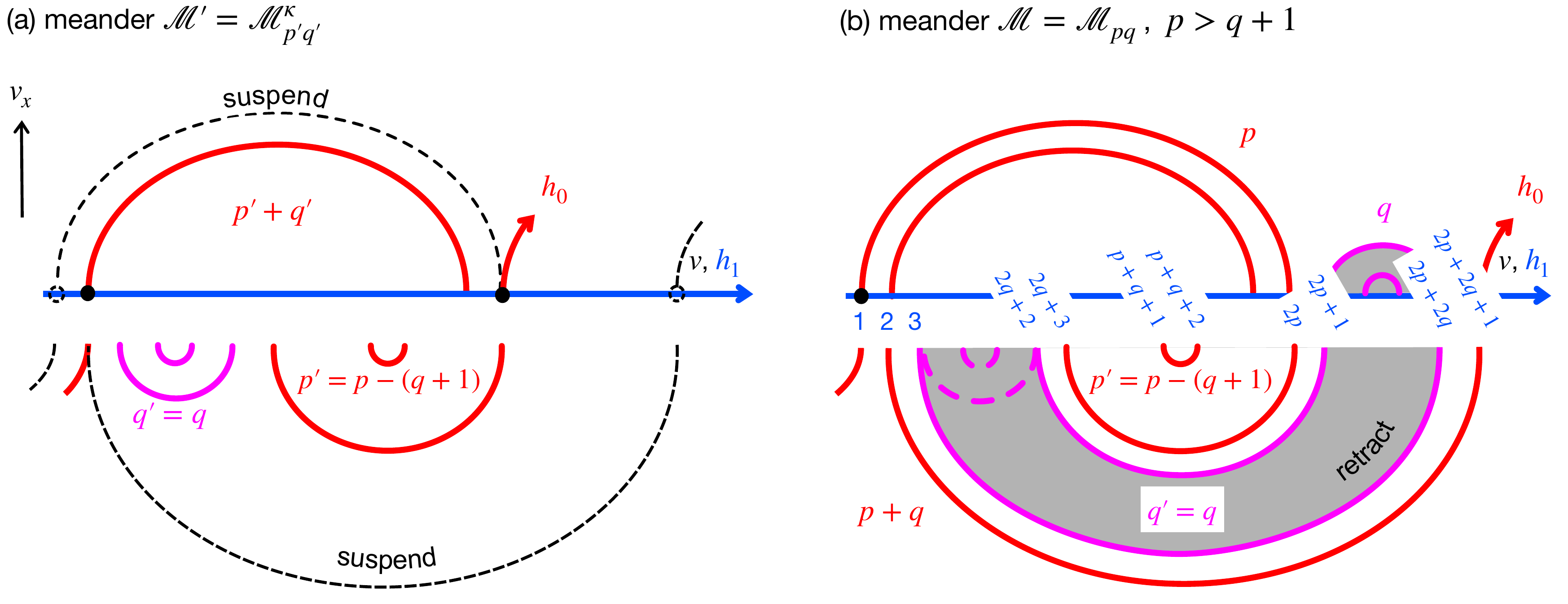}
\caption{\emph{
The case $r(q+1)<p<(r+1)(q+1)$: recursion from $\mathcal{M}_{p,q}$\,, (b), to the rotated meander $\mathcal{M}_{p',q'}^\kappa$\,, (a). 
As in \eqref{p'q'}, we consider $p'=p-(q+1), \  q'=q$ . 
(a) Suspension (black dashed) of $\mathcal{M}_{p',q'}^\kappa$ (red and purple) raises all Morse indices by 1; see proposition \ref{suspprop}(iv).
(b) Retraction of the $q$-nest (purple, shaded), from right to left, also raises all Morse indices by 1.
Read from left to right, conversely, the raising of Morse indices in the $q'$-nest, by suspension of (a), is compensated for, by the lowering left turn towards the $q$-nest, in (b).
This makes the Morse indices of corresponding equilibria in the $q$-nest of (b) and the $q'$-nest (a) coincide.
Reduction to the negative Morse index $i(A')=-1$, encountered after $r$ steps in the $q'$-nest, for $0<p'<q'+1$, provides a contradiction to the Morse property of $\mathcal{M}_{p,q}$\,.
}}
\label{qnose}
\end{center}
\end{figure}

To settle the case $r(q+1)<p<(r+1)(q+1)$ for $r\geq 1$, we proceed recursively; compare figs.~\ref{qnose}, (a) and (b).
In (b), consider the shaded area extending towards the upper $q$-nest of the meander $\mathcal{M}_{p,q}$\,.
Evidently, the retraction of the shaded area to the left produces a standard suspension of the meander $\mathcal{M}'$ in (a).
We also note $\mathcal{M}' = \mathcal{M}_{p',q'}^\kappa$ is trivially equivalent, by the rotation $\kappa$ of \eqref{rot}, to the meander $\mathcal{M}_{p',q'}$ with
\begin{equation}\label{p'q'}
p'=p-(q+1), \qquad q'=q.
\end{equation}
Proceeding from (a) to (b), now, suspension first raises every Morse number by 1; see proposition \ref{suspprop}(iv).
Next, the reinsertion of the lower left $q'$-nest of (a), by left turns, as the upper right $q$-nest of (b), reduces the Morse numbers of vertices in the reinserted $q$-nest by 1, compared to their originals in the $q'$-nest; see \eqref{i}.
In total, i.e. after suspension and reinsertion, the lower left $q'$-nest of $\mathcal{M}_{p',q'}$ and the upper right $q$-nest of $\mathcal{M}_{p,q}$ feature the same Morse indices.
Inductively, $r$ retraction steps reduce our meander to the previous case, where $0<p<q+1$.
Since the $q$-nest there did contain some $B$ or $B'$ with negative Morse index $i(B')=-1$, our recursive proof of theorem \ref{nonmorse}(ii) is now also complete.

\section{The simplest interesting case}\label{r=1}

In this section we address the remaining four theorems \ref{krmuthm}, \ref{krrhothm}, \ref{krballthm}, and \ref{krrevthm}, of section \ref{Res}, on the primitive 3-nose Sturm attractors $\mathcal{A}_{rq}$\,, their dissipative Morse meanders $\mathcal{M}_{r(q+1),q}$\,, and their entourage of Sturm permutations $\sigma_{rq}$ and connection graphs $\mathcal{C}_{rq}$\,.
For brevity and simplicity, we restrict our proofs to the simplest interesting case $r=1$.
We skip the trivial case $r=q=1$, already treated in fig.~\ref{susp}(a)-(c).
In subsection \ref{krrhothmpf} we use conspicuous nose locations to identify the action of trivial equivalences among these objects.
In particular we prove the trivial equivalence of $\mathcal{A}_{q1}$ and $\mathcal{A}_{1q}$ claimed in theorem \ref{krrhothm}, for $r=1$.
Theorem \ref{connthm} in subsection \ref{C1q} identifies the connection graphs.
This will easily prove the remaining three theorems, in subsections \ref{krmupf}--\ref{krrevpf}.
As an afterthought, we conclude with explicit expressions for the Sturm permutations $\sigma_{1q}$ and their trivially equivalent relatives, in proposition \ref{sigmaprop} of subsection \ref{sigmas}.

\subsection{Proof of theorem \ref{krrhothm}.}\label{krrhothmpf}

To locate noses $AB$ of equilibria we use the matrix notation $(a_1,b_1\,|\,a_0,b_0)$ for locations $a_\iota := h_\iota^{-1}(A),\ b_\iota := h_\iota^{-1}(B)$ and $\iota=0,1$.
Note how noses are characterized by adjacency $|a_\iota - b_\iota|=1$ under both boundary orders $h_\iota$\,.

\begin{lem}\label{noseloclem}
The following are corresponding nose locations of the indicated Sturm permutations, for any fixed $q\geq 2$:
\begin{enumerate}[(a)]
  \item the upper right nose $(4q+1,4q+2\,|\,2q+1,2q+2)$ of $\sigma_{q1}$;
  \item the lower left nose $(2,3\,|\,2q+2,2q+3)$ of $\sigma_{q1}^\kappa$;
  \item the nose $(2q+1,2q+2\,|\,4q+1,4q+2)$ of the upper rainbow of $\sigma_{1q}^\kappa$;
  \item the nose $(2q+2,2q+3\,|\,2,3)$ of the lower rainbow of $\sigma_{1q}$.
\end{enumerate}
The correspondence is under the trivial equivalences \eqref{rotsig}, for $\kappa$, and \eqref{invsig}, for $\rho$, as illustrated in fig.~\ref{equi}(a)-(d) for the special case $q=2$.
In particular, the four permutations are trivially equivalent and \eqref{krrho} holds, for $r=1$.
\end{lem}

\begin{proof}
The lower rainbow nose (d) of $\sigma_{1q}$, i.e.~for $p=q+1$, is obviously located at $(2p,2p+1\,|\,2,3) = (2q+2,2q+3\,|\,2,3)$, by arc counting.
Similarly, the upper rainbow nose (c) for the rotated meander associated to $\sigma_{1q}^\kappa$ is just as obviously located at the rotated position $(2p-1,2p\,|\,2p+2q-1,2p+2q) = (2q+1,2q+2\,|\,4q+1,4q+2)$.

Inversion $\rho$ of $\sigma$ interchanges the roles of $h_0$ and $h_1$. 
This swaps the entries of the nose matrix before and after the separator ``$|$''.
Therefore the noses corresponding to the rainbow noses in (c) and (d) become $(4q+1,4q+2\,|\,2q+1,2q+2)$ and $(2,3\,|\,2q+2,2q+3)$ in (a) and (b), respectively.
The first $h_1$-entries locate these noses at the extreme right and left of the horizontal $h_1$ axis, respectively.

It remains to show that the permutation $\sigma_{q1}$ in (a) is indeed the inverse of the Sturm permutation $\sigma_{1q}^\kappa$ in (c). 
(The other pair (b), (d) is treated analogously.)
From section \ref{Rot}, we already know that inversion $\rho$ preserves the number of noses and, up to $\kappa$, commutes with suspension; see \eqref{suspr}.
Therefore the inverse $\sigma_{1q}^{\rho\kappa}$ of $\sigma_{1q}^\kappa$ must also be a primitive 3-nose Sturm permutation $\sigma_{r'q'}$.
The upper nose in (a) is located rightmost, at $h_1\in \{4q+1,4q+2\}$, and hence cannot sit inside any larger nest.
Therefore $\sigma_{1q}^{\kappa \rho}=\sigma_{r'1}$ for some $r'$.
This implies $r'=q$, since the total number $2(p+q)+1=4q+3=4r'+3$ of vertices is also preserved under inversion $\rho$.
This proves the lemma, \eqref{krrho}, and theorem \ref{krrhothm}.
\end{proof}

\begin{figure}[t]
\begin{center}
\centering \includegraphics[width=\textwidth]{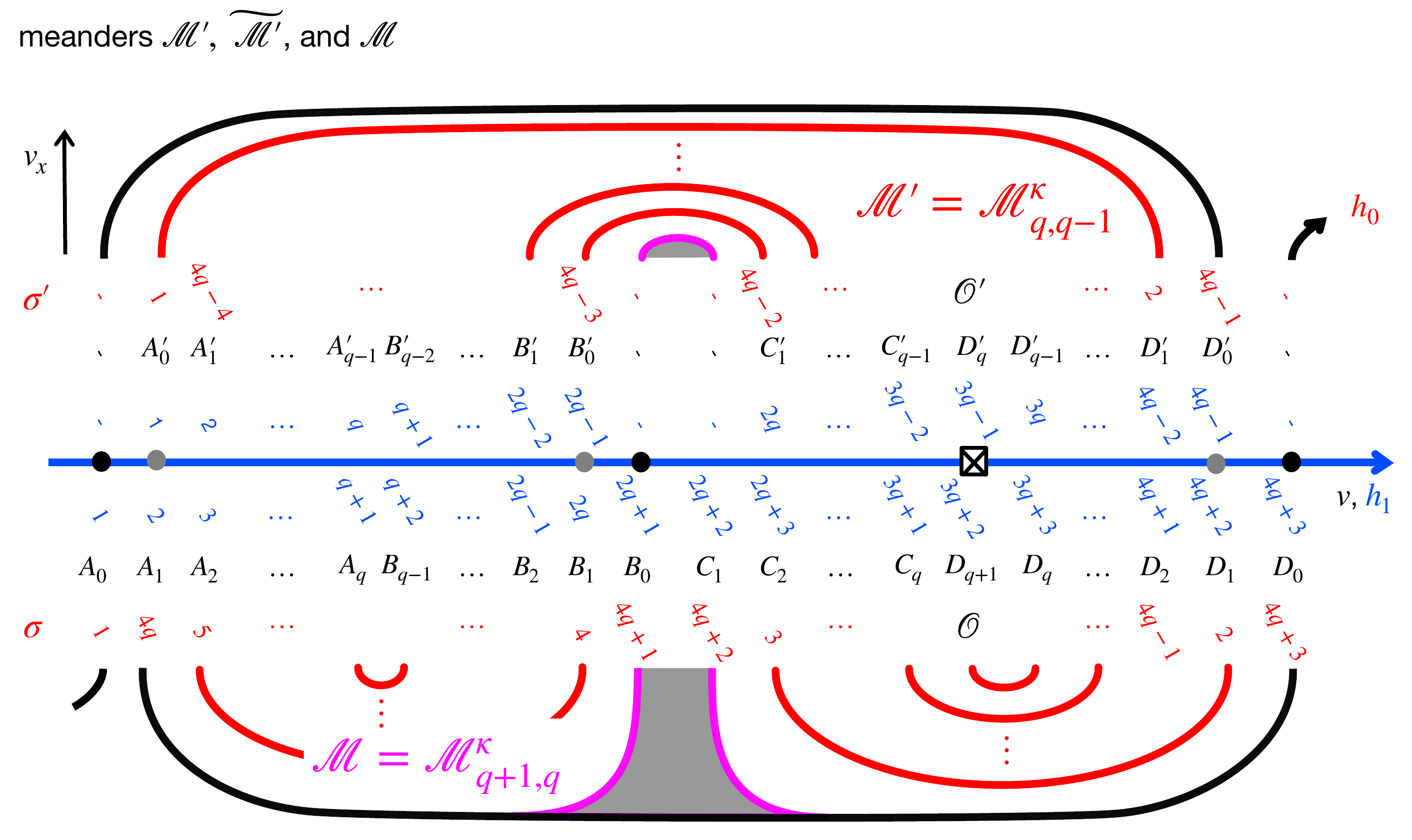}
\caption{\emph{The induction step $2\leq q-1 \mapsto q$: 
the meander $\mathcal{M}'=\mathcal{M}_{q,q-1}^\kappa$\,, its suspension $\widetilde{\mathcal{M}'}$, and the $q$-successor $\mathcal{M}=\mathcal{M}_{q+1,q}^\kappa$\,.
Below the horizontal $h_1$-axis are the $h_1$-order $j$ (blue), the equilibrium labels $A_0, \ldots$ of \eqref{connlab} (black), and the $h_0$-order $\sigma(j)$ (red), all for $\mathcal{M}$. 
Above the axis are the corresponding entries for $\mathcal{M}'$.
The two suspension arcs $A_0D_0'$ and $A_0'D_0$, from $\mathcal{M}'$ to $\widetilde{\mathcal{M}'}$, are indicated in black.
The nose insertion $B_0C_1$, from $\widetilde{\mathcal{M}'}$ to $\mathcal{M}$, is delineated by the purple boundary of the shaded region.
}}
\label{connmeander}
\end{center} 
\end{figure}

\subsection{The connection graphs $\mathcal{C}_{1q}$ and $\mathcal{C}_{q1}$}\label{C1q}

Next, we determine the connection graphs $\mathcal{C}_{1q}$ by recursion on $q$.
Since trivial equivalences induce isomorphism of connection graphs, and to simplify notation, we choose to compare the rotated connection graphs
\begin{equation}\label{connrec}
\mathcal{C}' := \mathcal{C}_{1,q-1}^\kappa \qquad \textrm{and} \qquad \mathcal{C} := \mathcal{C}_{1,q}^\kappa\,,
\end{equation} 
instead, along with their entourages of Sturm permutations $\sigma',\sigma$, meanders $\mathcal{M}', \mathcal{M}$, and attractors $\mathcal{A}',\mathcal{A}$.
See fig.~\ref{3ball}(a),(c),(d) for $\sigma, \mathcal{A}, \mathcal{C}$ in the simplest case $q=2$.
For general $q\geq 2$, we label the equilibria $h_1(j)$ of $\mathcal{A}$, etc., along the horizontal $h_1$-axis, as follows:
\begin{equation}\label{connlab}
E_j = h_1(j) =:
\begin{cases}
      A_{j-1}\,, & \textrm{for}\quad \qquad\,1 \leq j \leq q+1, \\
      B_{(2q+1)-j}\,, & \textrm{for}\quad\;\, q+2 \leq j \leq  2q+1, \\
      C_{j-(2q+1)}\,, & \textrm{for}\quad 2q+2 \leq  j \leq 3q+1, \\
      D_{(4q+3)-j}\,, & \textrm{for}\quad 3q+2 \leq  j \leq 4q+3.
\end{cases}
\end{equation}
Occasionally, we will also use the notation $\mathcal{O}:=D_{q+1}$\,.
Along the $h_1$-axis, this enumerates the equilibrium sequence by alternatingly ascending and descending subscripts as
\begin{equation}\label{h1equi}
A_0 \ldots A_q\; B_{q-1}\ldots B_0\; C_1 \ldots C_q\ (\mathcal{O}\textrm{=}D_{q+1})\ D_q\ldots D_0 \,.
\end{equation}
For $\mathcal{A}',\mathcal{C}'$ we use the corresponding notation $A_0'\,, \ldots$ . 
See figs.~\ref{connmeander}, \ref{conngraph} for illustration.


\begin{thm}\label{connthm}
Let $q\geq 2$. In the notation \eqref{connrec}-\eqref{h1equi}, the connection graph $\mathcal{C}$ is then given by
\begin{align}
\label{ato}
A_j \quad &\leadsto \quad \{A_{j-1},B_{j-1}\},\quad \textrm{for } 1\leq j \leq q;\\
\label{bto}
B_j \quad &\leadsto \quad \{A_{j-1},B_{j-1}\},\quad \textrm{for } 1\leq j \leq q-1;\\
\label{cto}
C_j \quad &\leadsto \quad \{B_{j-1},C_{j-1},D_{j-1}\},\quad \textrm{for } 2\leq j \leq q;\\
\label{dto}
D_j \quad &\leadsto \quad \{A_{j-1},C_{j-1},D_{j-1}\},\quad \textrm{for } 2\leq j \leq q+1;\\
\label{c1to}
C_1 \quad &\leadsto \quad \{B_0,D_0\};\\
\label{d1to}
D_1 \quad &\leadsto \quad \{A_0,D_0\}.
\end{align}
In particular, all admissible subscripts $j$ indicate Morse indices:
\begin{equation}\label{ilj}
i(L_j) = j,\quad \textrm{for all equilibrium tags } L\in\{A,B,C,D\}.
\end{equation}
\end{thm}

\begin{proof}
With the case $q=2$ already settled, we proceed by induction on $2\leq q-1 \mapsto q$.
We may therefore assume that the theorem already holds true for the ($q$-1) meander $\mathcal{M}'$ and its connection graph $\mathcal{C}'$, as illustrated in  figs.~\ref{connmeander} and \ref{conngraph}(a).
Starting from $q-1$, our first step is by suspension to $\widetilde{\mathcal{M}'}, \widetilde{\mathcal{C}'}$ as in figs.~\ref{connmeander} and \ref{conngraph}(b). Our second step, leading to the $q$-meander $\mathcal{M}$ and its connection graph $\mathcal{C}$, is by nose insertion; see figs.~\ref{connmeander} and \ref{conngraph}(c).

Suspension, our first step, invokes proposition \ref{suspprop}.
The equilibria $\widetilde{E}_j'$ of the suspension have been labeled $\widetilde{L}_j'$, to correspond to our notation $L_j'$ for $E_j'$\,.
Suspension raises Morse indices by 1, due to proposition \ref{suspprop}(iv).
Only for the cone vertices $\widetilde{E}_0'$ and $\widetilde{E}_{N+1}'$ of the suspension, at the lowest Morse level $i=0$, we have substituted the new labels $A_0,D_0$ in figs.~\ref{connmeander} and \ref{conngraph}(b).
The connection graph $\widetilde{\mathcal{C}'}$ of (b) then follows from the suspension corollary \ref{suspcor}.

\begin{figure}[t]
\begin{center}
\centering \includegraphics[width=0.97\textwidth]{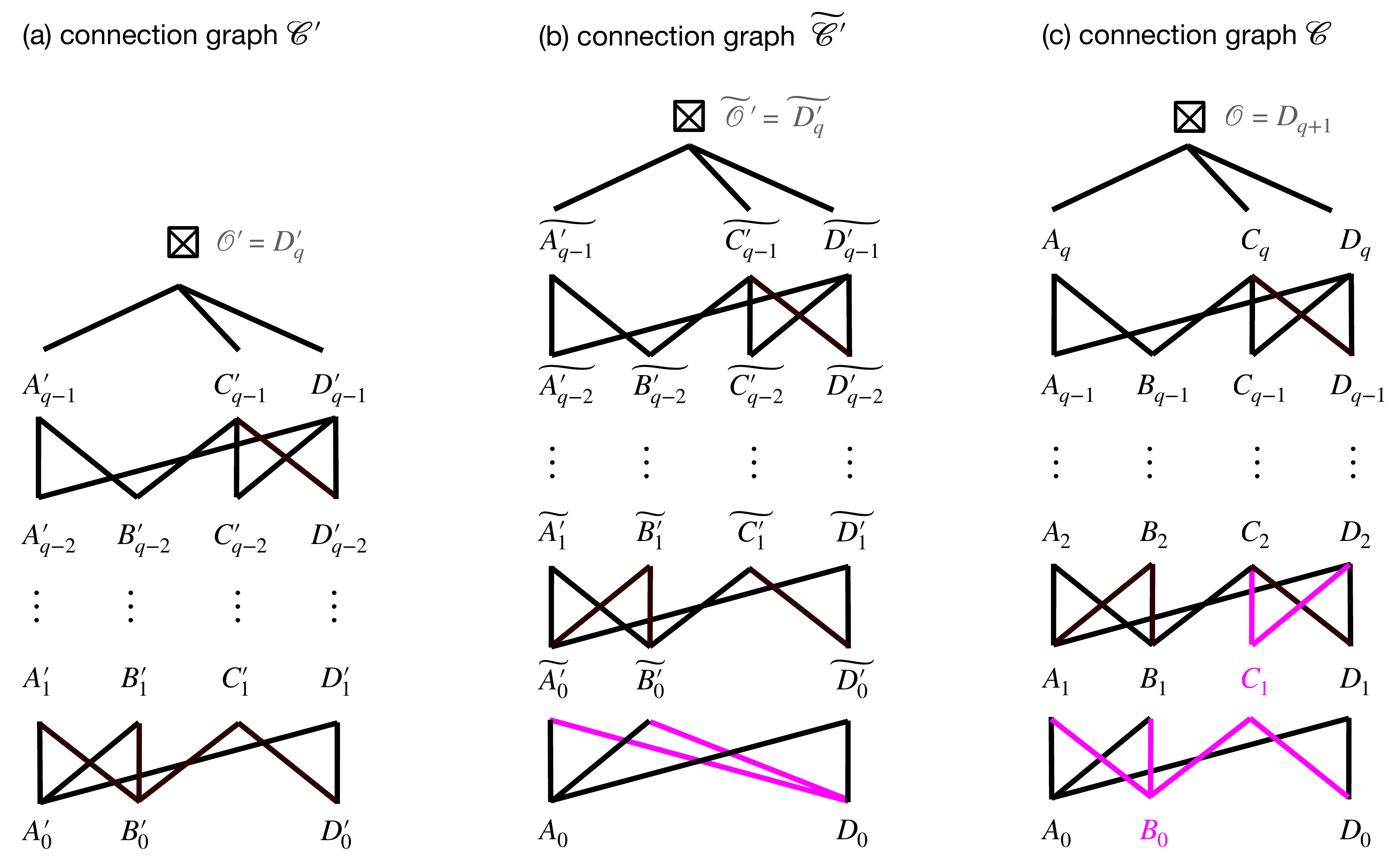}
\caption{\emph{The induction step $2\leq q-1 \mapsto q$: 
the connection graphs $\mathcal{C}'=\mathcal{C}_{1,q-1}^\kappa\,,\ \widetilde{\mathcal{C}'}$, and $\mathcal{C}=\mathcal{C}_{1,q}^\kappa$\,.
As before, adjacent rows indicate adjacent Morse levels.
Downwards arrows have been omitted.
(a) The connection graph $\mathcal{C}'$ of $\mathcal{M}'$, for $q-1$. Note how subscripts of equilibrium labels coincide with their Morse indices, from the lowest row sinks, at $i=0$, up to $i=q$ at the top unstable equilibrium $\mathcal{O}'=D_q'$\,.
(b) The double cone suspension from $\mathcal{C}'$ to $\widetilde{\mathcal{C}'}$ raises all Morse levels by 1.
Note the new $i=0$ sink vertices $A_0,D_0$, i.e. the vertices of the two suspension cones. Necessarily, both vertices are heteroclinic targets from any other equilibrium. 
In $\widetilde{\mathcal{C}'}$ this is manifest by the six heteroclinic arrows $\widetilde{A}_0', \widetilde{B}_0', \widetilde{D}_0' \leadsto \{A_0,D_0\}$.
(c) The connection graph $\mathcal{C}$ of $\mathcal{M}$, for $q$. 
Equilibria $\widetilde{L}_j'$ in $\widetilde{\mathcal{C}'}$ have been relabeled as $L_{j+1}$, for all tags $L\in\{A,B,C,D\}$, to sync subscripts $j$ with the raised Morse indices $i$.
Note how all heteroclinic edges are preserved at Morse levels $i\geq 1$, upon insertion of the nose arc $B_0C_1$.
The new heteroclinic edges of $\mathcal{C}$ to, and from, the purple nose $B_0C_1$ are indicated in purple.
The two heteroclinic edges of $\widetilde{\mathcal{C}'}$ which have been blocked in $\mathcal{C}$, by $B_0C_1$, have also been marked purple, in (b).
}}
\label{conngraph}
\end{center} 
\end{figure}

Our second step is the nose insertion of figs.~\ref{connmeander} and \ref{conngraph}(c).
First note our substitution $\widetilde{L}_j' \mapsto L_{j+1}$, for equilibria inherited by (c) from (b).
This ensures $i(L_j)=j$, for $j\geq 1$, as claimed in \eqref{ilj}.
The cone vertices $A_0,D_0$ have not been relabeled.
However, we now have to address three possible effects of the newly inserted nose arc $B_0C_1$ on heteroclinic edges (purple) in fig.~\ref{conngraph}:
\begin{enumerate}[(i)]
  \item previous edges of (b) blocked by nose equilibria $B_0, C_1$;
  \item new edges in (c) emanating from the nose $B_0, C_1$;
  \item new edges in (c) terminating at the nose $B_0, C_1$.
\end{enumerate}
Note $i(B_1)=1$, by suspension \eqref{ilj}, and $i(B_0)=0$ by \eqref{i}.

We start with blocking of type (i).
By \eqref{wolfrum},\eqref{block}, new blockings of $v_1 \leadsto v_2$, i.e. purple edges in (b), only arise through nose equilibria $w\in\{B_0,C_1\}$ which are located between other $v_1$ and $v_2$ along the meander order $h_0$ of $\mathcal{M}$, and which satisfy \eqref{block}.
Since $B_0C_1$ is a nose arc, blocking by $C_1=h_0(4q+2)$ is equivalent to blocking by $B_0=h_0(4q+1)$\,.
Except for the last equilibrium $D_0=h_0(4q+3)$, all equilibria $v=h_0(j)$ inherited by (c) from (b) have $h_0$-position $j$ less than the second to last $h_0$-position $4q+2$ of $C_1$\,.
Therefore, $C_1$ (or $B_0$, equivalently) cannot block any of the heteroclinic edges inherited from (b), by (c), except possibly for edges from $\{A_1,B_1,D_1\}$ to $D_0$\,.
The edge between $\widetilde{D}_0'=D_1$ and the sink $A_0$, for example, cannot be blocked, because they are $h_0$-neighbors on the suspension arc $A_0D_1$.
Similarly, the edge between the $h_1$-neighbors $\widetilde{D}_0'=D_1$ and the sink $D_0$ remains non-blocked.
However, $z(A_1-B_0)=z(B_0-D_0)=z(A_1-D_0)=0$ implies that $w=B_0$ blocks $\widetilde{A}_0'=A_1=:v_1\leadsto v_2:=D_0$\,.
Here and below we refer to \eqref{z} along the orders of $h_0$ or $h_1$, equivalently, for the calculation of zero numbers.
Similarly, $\widetilde{B}_0'=B_1=:v_1\leadsto v_2:=D_0$ is blocked by $B_0$ at $z=0$.
This settles the effects of blocking, (i).

Next, we address new heteroclinic edges (ii) emanating  from the nose.
Obviously, edges cannot emanate from the sink $i(B_0)=0$.
Just as obviously, $C_1$ connects heteroclinically to its nose neighbor $B_0$, and to its $h_0$-neighbor $D_0$.
However $z(C_1-B_0)=z(B_0-A_0)=z(C_1-A_0)=0$ implies that $w=B_0$ blocks $C_1=:v_1\leadsto v_2:=A_0$\,.
This identifies all edges emanating from the nose, (ii).

It only remains to address new heteroclinic edges (iii) terminating at the nose.
Consider the target $B_0$, first.
Obviously, there are heteroclinic edges towards the sink $i(B_0)=0$ from its $h_\iota$-neighbors $A_1,B_1,C_1$, all at Morse level $i=1$.
The hypothetical edge $D_1 \leadsto B_0$ is blocked by $C_2$, at $z=1$.
This settles the three edges towards target $B_0$.

Finally, consider the target $C_1$ of (iii).
We proceed by checking the potential sources $A_2, B_2, C_2$, $D_2$, in alphabetical order.
The hypothetical edge $A_2 \leadsto C_1$ is blocked by $B_0$, at $z=0$.
Indeed $z(A_2-D_0)=0$ implies $z(A_2-B_0)=z(B_0-C_1)=z(A_2-C_1)=0$.
Similarly, $B_0$ blocks the hypothetical edge $B_2 \leadsto C_1$, at $z=0$.
Obviously, there is a heteroclinic edge towards the saddle $i(C_1)=1$ from its $h_1$-neighbor $C_2$ at Morse level $i=2$.
To show $h_0(4q-1)=D_2 \leadsto C_1=h_0(4q+2)$, just note that the only equilibria $h_0$-between $D_2$ and $C_1$ are $h_0(4q)=A_1$ and $h_0(4q+1)=B_0$.
However, the latter pair precedes the former, along $h_1$, and therefore cannot be blocking.

This establishes the connection graph $\mathcal{C}$ of $\mathcal{M}$, as illustrated in figs.~\ref{connmeander} and \ref{conngraph}(c). 
By induction on $q$, the theorem is now proved.
\end{proof}

We can now prove the remaining three main theorems \ref{krmuthm},  \ref{krballthm}, and \ref{krrevthm}, for $r=1$.
We repeat that lemma \ref{noseloclem}, which already established theorem \ref{krrhothm}, allows us to base our proofs on the trivially equivalent connection graphs $\mathcal{C}=\mathcal{C}_{1q}^\kappa$\,, instead of $\mathcal{C}_{1q}$ itself.
All three theorems will become easy corollaries of theorem \ref{connthm}; see also fig.~\ref{conngraph}.
We conclude with an explicit proof of equivalence theorem \ref{krrhothm} which is independent of our more abstract approach via lemma \ref{noseloclem}. 
Instead, it will be based on a direct, explicit, and elementary computation of the Sturm permutations $\sigma_{q1}, \sigma_{1q}$,  $\sigma_{q1}^{-1}$, and $\sigma_{1q}^{-1}$, in proposition \ref{sigmaprop}.

\subsection{Proof of theorem \ref{krmuthm}.}\label{krmupf}
For any $q\geq 2$, the connection graph of theorem \ref{connthm} establishes the Morse counts 
\begin{equation}\label{1qmu}
m_{1q}(i)\ =\ 
\begin{cases}
      3,& \textrm{for}\quad i=0;\\
      4,& \textrm{for}\quad 1\leq i \leq q-1;\\
      2(q+1-i)+1,& \textrm{for}\quad q\leq i\leq q+1.
\end{cases}
\end{equation}
See also fig.~\ref{conngraph}(c), and \eqref{oto} below, more explicitly.
For $r=1$, this proves the Morse counts $m_{rq}(i)$ of theorem \ref{krmuthm}. \hfill $\bowtie$

\subsection{Proof of theorem \ref{krballthm}.}\label{krballpf}
By the Schoenflies theorem \cite{firo15}, it is sufficient to prove that the single equilibrium $\mathcal{O}=D_{q+1}$ of the top Morse index $i(\mathcal{O})=q+1=\dim \mathcal{A}_{1q}$ connects heteroclinically to all other equilibria $E$. 
In symbols, $D_{q+1}\leadsto E$.
By transitivity of the directed edge relation $\leadsto$, this amounts to showing that there exists a di-path from $D_{q+1}$ to any $E$, in the connection di-graph $\mathcal{C}$.
This is obvious from \eqref{ato}--\eqref{d1to}, which coarsen to
\begin{equation}\label{oto}
\begin{aligned}
\mathcal{O}=D_{q+1}\ &\leadsto\ \{A_q,C_q,D_q\}\ \leadsto\ \{A_{q-1},B_{q-1},C_{q-1},D_{q-1}\}\ \leadsto \\
&\leadsto\ \ldots\ \leadsto\ \{A_1,B_1,C_1,D_1\}\ \leadsto\ \{A_0,B_0,D_0\}\,. 
\end{aligned}
\end{equation}
Indeed, all equilibria $E$ from \eqref{h1equi} occur in this sequence. \hfill $\bowtie$

\subsection{Proof of theorem \ref{krrevthm}.}\label{krrevpf}
For the connection graph $\mathcal{C}$ of theorem \ref{connthm}, fig.~\ref{conngraph}(c), consider the involutive vertex map
\begin{equation}\label{krrev}
\mathcal{R}:\ \quad A_j\, \longleftrightarrow\, D_{q-j}\,,\ \quad B_j\, \longleftrightarrow\, C_{q-j}\,.
\end{equation}
Here $0\leq j \leq q$ in the first swap, but $0\leq j \leq q-1$ in the second.
Inspection shows that $\mathcal{R}$ is a reversor automorphism of the connection di-graph $\mathcal{C}\setminus D_{q+1}\,$.
Indeed, $\mathcal{R}$ reverses all arrows.  \hfill $\bowtie$

\subsection{Explicit Sturm permutations}\label{sigmas}

We derive the explicit primitive 3-nose Sturm permutations $\sigma_{q1},\sigma_{q1}^{-1},\sigma_{1q},\sigma_{1q}^{-1} \in S_{4q+3}$\,.

\begin{prop}\label{sigmaprop}
Claim \eqref{krrho} of theorem \ref{krrhothm} holds true, for $r=1$ and all $q\geq1$, due to the following explicit expressions of the relevant permutations.
\begin{enumerate}[(i)]
  \item With arguments $0\leq j \leq q$, as appropriate, the permutation $\sigma_{q1}$ satisfies
\begin{equation} \label{sigq1-0}
  \begin{aligned}
  \sigma_{q1}(4j)\quad\ \; &= 2q-2j+2;\\
  \sigma_{q1}(4j+1) &= \qquad\, 2j+1;\\
  \sigma_{q1}(4j+2) &= 4q-2j+2;\\
  \sigma_{q1}(4j+3) &= 2q+2j+3.
  \end{aligned}
\end{equation}
  \item The inverse permutation $\sigma_{q1}^{-1}$ is given explicitly by
    \begin{align}
   \label{sigq1invodd}
 \sigma_{q1}^{-1}(2j'+1) &=
 \begin{cases}
      4j'+1,& \text{for } 0\leq j' \leq q; \\
      4(j'-q-1)+3,\phantom{i} & \text{for } q< j' \leq 2q+1;
\end{cases}
  \\
   \label{sigq1inveven}
 \sigma_{q1}^{-1}(2j')\quad\ \; &=
 \begin{cases}
      4(q+1-j'),& \text{for } 1\leq j'\leq q; \\
      4(2q+1-j')+2,& \text{for } q<j' \leq 2q+1.
\end{cases}
  \end{align}
  \item The permutation $\sigma_{1q}$ is given explicitly by
    \begin{align}
   \label{sig1qodd}
 \sigma_{1q}(2j'+1) &=
 \begin{cases}
      4j'+1,& \text{for } \quad\ 0\leq j' \leq q; \\
      4(j'-q-1)+3,& \text{for } \quad\ q< j' \leq 2q+1;
\end{cases}
  \\
   \label{sig1qeven}
 \sigma_{1q}(2j')\quad\ \; &=
 \begin{cases}
      4(q+1-j')+2,& \text{for }\quad\,\  1\leq j'\leq q+1; \\
      4(2q+2-j'),& \text{for } q+1<j' \leq 2q+1.
\end{cases}
  \end{align}
  \item With arguments $0\leq j \leq q$, as appropriate, the permutation $\sigma_{1q}^{-1}$ satisfies
\begin{equation}\label{sig1qinv-3}
 \begin{aligned}
  \sigma_{1q}^{-1}(4j)\quad\ \; &= 4q-2j+4;\\
  \sigma_{1q}^{-1}(4j+1) &= \qquad\,2j+1;\\
  \sigma_{1q}^{-1}(4j+2) &= 2q-2j+2;\\
  \sigma_{1q}^{-1}(4j+3) &= 2q+2j+3.
  \end{aligned}
\end{equation}
\end{enumerate}
\end{prop}

\begin{proof}
Obviously, the sixteen expressions \eqref{sigq1-0}-\eqref{sig1qinv-3} define permutations in $S_{4q+3}$\,.
Just for the moment, let us denote by $\sigma_{q1}$ and $\sigma_{q1}^{-1}$ the expressions in (i),(ii).
Then $\sigma_{q1} \circ \sigma_{q1}^{-1} = \mathrm{id}$ is obvious, by definition. 
Therefore (ii) actually defines the inverse permutation of (i).

To show that $\sigma_{q1}^{-1}$ actually is the inverse permutation of the true meander permutation $\sigma_{q1}$ associated to the dissipative meander $\mathcal{M}_{2q,1}$, we first recall that, equivalently, $\sigma_{q1}^{-1}$ is supposed to provide the correct vertex locations, along the horizontal $h_1$-axis, if we enumerate vertices along the meander curve by $h_0= \mathrm{id}$.
Arcs of the meander take the form $h_0(j)\,h_0(j+1)$.
Since the meander switches sides, at each crossing with the horizontal axis, we obtain upper arcs for odd $j$. 
Lower arcs have even $j=2j',\ 1\leq j' \leq 2q+1$.
The explicit forms \eqref{sigq1inveven},\eqref{sigq1invodd} then imply the invariance
\begin{equation}\label{lowrain}
 \sigma_{q1}^{-1}(2j')+\sigma_{q1}^{-1}(2j'+1)=4q+5,
\end{equation}
for all $j'$. 
This characterizes the lower arcs of $\mathcal{M}_{2q,1}$ to form a nested rainbow.

For the two upper nests of $\mathcal{M}_{2q,1}$, we argue analogously for the upper arcs from odd $j=2j'+1,\ 0\leq j' \leq 2q$ to even $j+1=2(j'+1)$.
Note first how $j'=q$ defines the upper 1-nest of the rightmost nose at horizontal positions $4q+1, 4q+2$; see \eqref{sigq1invodd},\eqref{sigq1inveven}, and lemma \ref{noseloclem}(a).
For all remaining $j'\neq q$, we obtain the invariant
\begin{equation}\label{upnest}
 \sigma_{q1}^{-1}(2j'+1)+\sigma_{q1}^{-1}(2(j'+1))=4q+1,
\end{equation}
which characterizes the left upper 2q-nest of $\mathcal{M}_{2q,1}$.
This establishes our formal expression of $\sigma_{q1}$, in (i), to be the true meander permutation of the meander $\mathcal{M}_{2q,1}$, as introduced in the first paragraph of section \ref{Res}.
The argument for $\sigma_{1q}$, (iii), proceeds analogously, and is left as a useful exercise.

We are now able to check \eqref{krrho} for $r=1$, by brute force instead of our previous deeper insight.
We just evaluate the left hand side of the equivalent claim $(\kappa\sigma_{q1}^{-1} \kappa) (j)=\sigma_{1q}(j)$, separately, for even and odd arguments $j$.
We only present the straightforward calculation for the even case $j=2j',\ 1\leq j' \leq q+1$; the three other cases are analogous.
Inserting the flip $\kappa(j)=4q+4-j$ from \eqref{flip}, and definition \eqref{sigq1inveven} of $\sigma_{q1}^{-1}$, we obtain
\begin{equation}\label{1qrhocheck}
\scalebox{0.97}[1.0]{$
\begin{aligned}
(\kappa\sigma_{q1}^{-1} \kappa) (2j') &= (\kappa\sigma_{q1}^{-1})\big( 2(2q+2-j')\big) =\kappa \big( 4\big( 2q+1-(2q+2-j')\big) +2\big) = \\
&=\kappa (4j'-2) = 4q+4-(4j'-2)=4(q+1-j')+2 = \sigma_{1q}(2j'),
\end{aligned}
$}
\end{equation}
in agreement with definition \eqref{sig1qeven} of $\sigma_{1q}(2j')$, for $1\leq j' \leq q+1$.

Finally, we obtain (iv) via $\sigma_{1q}^{-1}=\kappa\sigma_{q1}\kappa$. Alternatively we may check the inversion (iv) of (iii) formally, as we did for the pair (i),(ii).

This proves the proposition.
\end{proof}

\section{Discussion}\label{Dis}

We discuss some broader settings for our results.
See subsection \ref{Dis-r} for the cases $r>1$ of our main results in section \ref{Res}, which section \ref{r=1} did not provide proofs for.
In \ref{Dis-gen} we briefly mention some pertinent literature on fully nonlinear equations, grow-up, and blow-up.
ODE variants of the PDE \eqref{PDE}, like cyclic monotone feedback systems and Jacobi systems, arise by finite difference discretization. See subsection \ref{Dis-jac}.
In \ref{Dis-rev}, we conclude with some more topological aspects of our results, and the open question of time reversal for full boundary spheres $\Sigma=\partial\mathcal{A}$ of global attractors, rather than for just their connection graphs $\mathcal{C}\setminus\mathcal{O}$.

\subsection{The cases $r>1$}\label{Dis-r}

The proof of theorem \ref{nonmorse} in section \ref{Nonmorse} gives an indication on how to proceed inductively for $r>1$; see fig.~\ref{qnose}.
Of course, we may perform $q$ successive nose insertions as in fig.~\ref{qnose}(b) for $p=(r+1)(q+1)$ as well, coming from $p'=r(q+1)$.
In case $q=1$, this inserts just one nose of two equilibria, reminiscent of -- but, technically, slightly different from -- our insertion of the nose $B_0C_1$ in figs.~\ref{connmeander} and \ref{conngraph}(c).
That insertion occurred at Morse levels $i=0,1$.
In case $q>1$, more ambitiously, we are inserting a $q$-nest of $2q$ equilibria, at the lowest Morse levels $i=0,\ldots,q$.
This makes it more demanding, technically and notationally, to perform the requisite induction step $r\mapsto r+1$ for the connection graphs $\mathcal{C}_{rq}$.
As our starting point $r=1$, for any $q>1$, however, we may use the connection graphs $\mathcal{C}=\mathcal{C}_{1q}^\kappa$ already established in theorem \ref{connthm} and fig.~\ref{conngraph}(c).
We postpone the details to our sequel \cite{firo23}.

\subsection{Nonlinear PDEs, grow-up, and blow-up}\label{Dis-gen}

Technical groundwork for generalizations to fully nonlinear equations, including nonlinear boundary conditions, has been laid by Lappicy and coworkers \cite{la18,la20,lafi18,la21,labe22}.
An interesting class of geometric applications are curve-shortening flows in the plane \cite{an91}.

The qualitative behavior of parabolic global ``attractors'' of non-dissipative nonlinearities $f$ is a very intriguing subject, even in the semilinear case.
For general \emph{blow-up} in finite time see the monograph \cite{quitsoup} and, in the Sturm setting, also \cite{ga04}.
For an attempt to describe the development of sign-changing blow-up profiles by zero numbers, in one space dimension, see \cite{fima07}.
Alternatively, solutions may exhibit \emph{grow-up} to infinity, taking infinite time.
The set of bounded global solutions $u(t,.),\ t\in \mathbb{R}$, of \eqref{PDE} will still consist of only equilibria and heteroclinic orbits.
The question how global solutions may connect to infinity, ``heteroclinically'', is attracting increasing attention; 
see for example \cite{ben,pim,lapim,carpim} and the references there.

\subsection{Jacobi systems}\label{Dis-jac}

\emph{Jacobi systems} are a spatially discrete analogue, including a zero number dropping property \eqref{zdrop}; see \cite{fuol88}.
Motivated by, but much more general than, a semi-discretized finite-difference version of the PDE \eqref{PDE}, they take the ODE form
\begin{equation}\label{jac}
\dot{u_j} = f_j(u_{j-1},u_j,u_{j+1}),
\end{equation}
for $j=1,\ldots,n$. 
The partial derivatives of $f_j$ with respect to the off-diagonal entries $u_{j\pm 1}$ are assumed strictly positive.
See \cite{fibh} for an application to strongly damped mechanical oscillators.
\emph{Cyclic monotone feedback systems} are a limiting case of spatially periodic subscripts $j$ mod $n$, with $f_j$ independent of $u_{j+1}$.
See for example \cite{smi,mpsmi,fi21}, with applications to gene feedback cycles and oscillations,
and \cite{mpse} for an extension which includes an additional time-delay.
For ``Neumann'' (or other separated) boundary conditions like $u_0:=u_1,\ u_{n+1}:=u_n$\,, but not necessarily for periodic boundary conditions, the system is still gradient-like \cite{fige}.

Equilibria of Jacobi systems satisfy the recursion relation $0=f_j(u_{j-1},u_j,u_{j+1})$.
Solving these equations for $u_{j+1}$, implicitly, equilibria satisfy an equivalent 2-term recursion
\begin{equation}\label{jac0}
u_{j+1} = g_j(u_j,u_{j-1}),
\end{equation}
for $j=1,\ldots,n$.
Here the partial derivative of $g_j$ with respect to $u_{j-1}$ becomes strictly negative.
In the dissipative case, gobal attractors $\mathcal{A}$, connection graphs $\mathcal{C}$, shooting meanders $\mathcal{M}$, and meander permutations $\sigma$ can be defined in complete analogy to the PDE case, with orbit equivalent attractors for equal permutations \cite{firo00}.
The role of the horizontal $h_1$-axis is then taken over by the diagonal $u_{n+1}=u_n$\,, in the $(u_n,u_{n+1})$-plane.

Concrete nonlinearities are still largely unexplored, even for the case where $g_j=g$ does not depend on $j$.
For a prominent example, we mention the Chirikov standard map, often written as
\begin{equation}\label{chirpu}
\begin{aligned}
p_{j+1} &= p_j + \lambda \sin u_j\,,\\
u_{j+1} &= u_j + p_{j+1}\,.
\end{aligned}
\end{equation}
Eliminating $p_j=u_j-u_{j-1}$ in \eqref{chirpu}, we obtain the equivalent 2-term recursion
\begin{equation}\label{chirrec}
u_{j+1} = g_j(u_j,u_{j-1}) := 2u_j +  \lambda \sin u_j - u_{j-1}\,.
\end{equation}
with negative partial derivative as required in \eqref{jac0}.

Similarly, the celebrated Hénon map is a 2-term recursion
\begin{equation}\label{henrec}
u_{j+1} = g_j(u_j,u_{j-1}) := 1-a u_j^2 + b u_{j-1}\,.
\end{equation}
Negative partial derivative, as in \eqref{jac0}, requires $b<0$.
For the usual sign, $b>0$, we may consider even $n$, and revert the signs of every other \emph{pair}, i.e. $u_j\mapsto -u_j$, for each $j\equiv 2,3 \mod 4$.
That recovers a properly negative partial derivative of $g$ with respect to $u_{j-1}$.
For odd length $n$, we can define meanders and permutations with respect to the ``anti-Neumann'' off-diagonal $u_{n+1}=-u_n$, in the $(u_n,u_{n+1})$-plane, rather than the diagonal.
The specific meander permutations which arise in such standard examples, by shooting, have never been addressed in any systematic way, to our knowledge.
For some related remarks in the context of Anosov maps see \cite{fi05}.

\subsection{Time reversal and reversibility}\label{Dis-rev}

One elementary formal operation on a Sturm meander $\mathcal{M}$ is a vertical flip, to some meander $\mathcal{W}$, by reflection at the horizontal $h_1$-axis. 
Let $\sigma,\varrho$ denote the associated meander permutations, respectively.
The flipped meander $\mathcal{W}$ emanates below the horizontal axis, from vertex $h_0(1)=h_1(1)$, but remains otherwise  dissipative, formally.
Inspection of Morse numbers \eqref{i}, however, now replaces any $i_j$ by $-i_j$\,.
Indeed, right turns on $\mathcal{M}$ become left turns, on $\mathcal{W}$.
Induced by $z_{kk}=i_k$ in \eqref{z}, the zero numbers also reverse sign.
Adjacency and blocking \eqref{block}, however, remain unaffected.
In terms of formal connectivity \eqref{wolfrum}, the sign reversal of the Morse numbers $i$ reverses all arrows in the associated formal connection graph of $\mathcal{W}$.
However, what does such algebraic trickery mean, in terms of actual equilibria $v$ of \eqref{PDE},\eqref{ODE}, which cannot possibly possess negative Morse indices $i(v)$?

We have already observed in section \ref{Rot}, how repeated suspensions raise Morse numbers and zero numbers, but preserve formal connectivity; see proposition \ref{suspprop} and corollary \ref{suspcor}.
Let $d=\dim \mathcal{A}= \max_j i_j$ denote the dimension of the original Sturm attractor associated to $\mathcal{M}$.
Then $\mathcal{W}$ has -$d$ as its minimal Morse number.
Therefore $\mathcal{W}$ becomes Sturm, first, after $s=d$ suspensions.
In fact, the suspended connection di-graph of $\mathcal{W}$ will contain the time reversed, i.e. inversely oriented, original connection di-graph $\mathcal{C}$, as a subgraph; see proposition \ref{suspprop}(vii).

Alas, such \emph{time reversal} does not provide \emph{time reversibility}, i.e.~an involutive time reversor $\mathcal{R}$ \emph{within} one and the same connection graph, as we have encountered on the boundary spheres of $\mathcal{C}_{rq}$ within the Sturm balls $\mathcal{A}_{rq}$\,.
Surprising as time reversal may be, it only shows how
\begin{description}
  \item[\phantom{max}] \emph{any time-reversed Sturm connection graph appears within some larger Sturm connection graph, of the same dimension $d$.}
\end{description}

Time reversibility of the connection graph on the boundary sphere $\Sigma^{d-1} = \partial \mathcal{A}$, however, is not all that exceptional either, for Sturm balls $\mathcal{A} = \mathrm{clos} W^u(\mathcal{O})$ of dimension $d$.
Before the 3-nose examples $\mathcal{A}_{rq}$\,, with $d=r+q$, we had already encountered the Chafee-Infante balls $\mathcal{A}_d$\,. 
Other examples of such $\mathcal{A}$ are all planar $n$-gons \cite{firo2d-2}, and the solid tetrahedron \cite{firo14}.
The methods of \cite{firo3d-3} provide any self-dual graph on $\Sigma^2$, and any solid $d$-simplex.

All these examples exhibit a weak reversor $\mathcal{R}$, as in \eqref{rev}, just on the vertices of $\Sigma^{d-1} = \partial \mathcal{A}$ in the connection graph $\mathcal{C}$.
However, such a weak reversor $\mathcal{R}$ need not extend to a strong reversor, automatically, on all of $\Sigma^{d-1}$.
Indeed, $\mathcal{R}$ only establishes that  certain heteroclinic orbits possess twins, someplace else, which run in reverse direction.

To prove strong time reversibility, the strong reversor $\mathcal{R}$ needs to define an orbit equivalence, mapping PDE orbits to PDE orbits, but reversing their time direction.
Poincaré (self-)duality of the Thom-Smale complex \eqref{S} may provide a first step in that direction.
Reversing time in fact interchanges the roles of stable and unstable manifolds, in the Thom-Smale complex.
Although this seems problematic in our infinite-dimensional PDE setting, finite-difference semi-discretization allows us to consider Jacobi systems \eqref{jac}, instead, where duality becomes finite-dimensional. 
Alternatively, we could work inside the global attractor $\mathcal{A}$ itself.
In general, duality on $\mathcal{A}$ is fraught with the additional difficulty that $\mathcal{A}$ may contain pieces of different local dimension.
The above cases of a sphere boundary $\partial\mathcal{A}_{rq}=\Sigma^{r+q-1}$ have been more benign.

There remain several steps towards elusive strong reversibility on $\Sigma^{r+q-1}$.
We have to show that the dual Thom-Smale complex is equivalent to its original, at least combinatorially.
We then have to establish topological equivalence of the complexes.
And, finally, we have to carefully adapt the duality construction such that the topological equivalence actually maps orbits to orbits, in the underlying Jacobi system.

Only in the special Chafee-Infante case of section \ref{ChIn} it seems fairly clear how to achieve that.
The general task certainly lies beyond the scope of the present paper.
Instead, we present a simple ODE model which, at least, features the same reversible connection graph $\mathcal{C} = \mathcal{C}_{1q}^\kappa$ of fig.~\ref{conngraph}(c), as the global attractor $\mathcal{A}=\mathcal{A}_{1q}^\kappa$ does, on the boundary sphere $\Sigma^q$.
An intriguing feature of $\mathcal{C}\setminus \mathcal{O}$ on $\Sigma^q$
are the two full Chafee-Infante sub-graphs $\mathcal{C}_q$\,: 
one with tags $L_j,\ L\in\{A,B\}$, on the left, and the other -- upside down, i.e. time-reversed -- on the right with tags $L \in \{C,D\}$.
Compare \eqref{CdO},\eqref{CdAB} and fig.~\ref{ci3}, with \eqref{ato}--\eqref{d1to} and fig.~\ref{conngraph}(c).

We therefore follow the Chafee-Infante model \eqref{cimi} of \cite{mi95}, to develop a model for the dynamics in $\Sigma^q$.
By stereographic projection, we represent $\Sigma^q$ as $\mathbb{R}^q$, with one-point compactification at infinity.
In $\mathbb{R}^q$, we choose polar coordinates $(\mathbf{r},\Phi)$.
In angular direction $\Phi$, as in \eqref{cimi}, we start from the diagonal matrix $\mathbf{Q} e_j=\mu_j e_j\,,\ j=0,\ldots, q-1$, with descending real eigenvalues $\mu_0>\ldots>\mu_{q-1}$\,.
As before, $\mathbf{P}_\Phi$ denotes orthogonal projection onto the tangent space at $\Phi$ in the unit sphere $S^{q-1}\subset \mathbb{R}^q$.
With the orthogonal involution $\mathbf{R} e_j := e_{q-1-j}$\,, our model then reads
\begin{equation}\label{revmodel}
\begin{aligned}
\dot{\mathbf{r}}    &= s(\mathbf{r})\, \mathbf{r}\,(\mathbf{r}^{-2}-4)(4-\mathbf{r}^2)\,;   \\
\dot{\Phi}&= s(\mathbf{r})\,\mathbf{P}_\Phi\left( (4-\mathbf{r}^2)\,\mathbf{Q}\Phi + (\mathbf{r}^{-2}-4)\,\mathbf{R}\mathbf{Q}\mathbf{R}\Phi \right)\,.
\end{aligned}
\end{equation}
We have scaled time by the multiplier $s(\mathbf{r}):=1/(\mathbf{r}^2+\mathbf{r}^{-2})$ for global existence, and to ensure regularity including $\mathbf{r}=0,\infty$, i.e. regularity on $\Sigma^q$.
Note symmetry $s(\mathbf{r})=s(1/\mathbf{r})$ and strong time reversibility, under the reversor 
\begin{equation}\label{revrev}
\mathcal{R}: \quad(\mathbf{r},\Phi) \mapsto (1/\mathbf{r},\mathbf{R}\Phi)\,.
\end{equation}
Indeed, $(\mathbf{r}(t),\Phi(t))$ solves \eqref{revmodel} if, and only if, $\mathcal{R}(\mathbf{r}(-t),\Phi(-t))$ does.

For $0\leq j < q$ we label the equilibria at positions $(\mathbf{r},\Phi)$ as follows:
\begin{equation}\label{revequi}
\begin{aligned}
    A_j &= (\tfrac{1}{2},e_j),      \qquad B_j = (\tfrac{1}{2},-e_j), \\
    C_{j+1} &= (2,-e_j),         \ D_{j+1} = (2,e_j), \\
        A_q &= \;0,           \quad  \ \ \,      \qquad D_0 = \;\infty\,. 
\end{aligned}
\end{equation}
The region $0\leq \mathbf{r}\leq\tfrac{1}{2}$ then models the standard Chafee-Infante sub-graph $\mathcal{C}^q$, with tags $L_j,\ L\in\{A,B\}$, analogously to \cite{mi95}.
By reversibility \eqref{revrev}, the time reversed flow in the region $2\leq \mathbf{r}\leq\infty$ models the time reversed, upside down, Chafee-Infante sub-graph with tags $L_j,\ L\in\{C,D\}$.
Indeed, subscripts reflect Morse indices in the original $q$-sphere $\Sigma^q$.
For example, the unstable eigenspaces of equilibria $A_j,B_j,C_{j+1},D_{j+1}$, in the angular directions $\Phi$, are all spanned by the preceding unit vectors $e_0,\ldots, e_{j-1}$.
In $\Sigma^q$ this amounts to Morse-Smale transversality of their stable and unstable manifolds.
The radial flow on each invariant eigenspace $\Phi = e_j$ provides the remaining heteroclinic orbits $C_{j+1} \leadsto B_j$ and $D_{j+1} \leadsto A_j$ of the connection graph $\mathcal{C}_{1q}^\kappa$ from fig.~\ref{conngraph}(c).
In our geometric model, at least, this glues the Chafee-Infante part $\mathcal{C}_q$\,, with tags $L_j,\ L\in\{A,B\}$, to the time-reversed Chafee-Infante part $\mathcal{C}_q$ with tags $L \in \{C,D\}$.

The basins of attraction are easily described, for the three $i=0$ sinks $A_0,B_0,D_0$\,. 
We only describe the basins within the sphere $\Sigma^q$, stereographically projected to $\mathbb{R}^d$.
The basin boundary of $D_0=\infty$ is the invariant ($q$-1)-sphere $\mathbf{r}=2$.
The invariant hyperplane $\langle e_1,\ldots,e_{q-1}\rangle$, inside the $q$-ball of radius $\mathbf{r}=2$, is the shared basin boundary of the other two sinks $A_0,B_0$, at angles $\Phi=\pm e_0$.
The same hyperplane splits the $\mathbf{r}=2$ sphere into two closed hemispheres, which are the shared basin boundaries of $A_0,B_0$ with $D_0$, respectively.
The equilibria in their intersection, i.e.~in the equatorial  $\mathbf{r}=2$ sphere of dimension $q-2$, and only those, possess heteroclinic orbits to all three sinks.

It is a useful exercise to locate all those equatorial and hemisphere equilibria, in our geometric desciption, and to verify their heteroclinic orbits to the respective sinks in the connection graph of fig.~\ref{conngraph}(c).
Indeed, all other equilibria with tags $C,D$ connect to $D_0$. 
Similarly, all non-sink equilibria of any tags, except $C_1$, connect to $A_0$\,. 
For $B_0$ the analogous exception is $D_1$\,.

Upon time reversal, equilibria in basin boundaries of sinks become heteroclinic targets of sources, instead. Therefore we can read off the basin boundary equilibria of $D_0,A_0,B_0$ from the targets of their reversor-related sources $A_q,D_q,C_q$ in fig.~\ref{conngraph}(c), respectively. 
See \eqref{revrev} and \eqref{revequi} for the precise tags and subscripts involved.
Again we see how, analogously, all other equilibria with tags $A,B$ are targets of $A_q$. 
Similarly, all non-source equilibria of any tags, except $B_{q-1}$, are targets of $D_q$\,. 
For $C_q$ the analogous exception is $A_{q-1}$\,.

Whether or not the same detailed geometry describes strong time reversibility in the sphere boundary of the Sturm attractor, $\partial\mathcal{A}_{1q}^\kappa = \Sigma^q$, remains open at present.
The explicit reversor $\mathcal{R}$ of \eqref{revrev} in our model \eqref{revmodel} certainly cautions us that this is not a trivial task.
But this is just one of the many curiosities, intricacies, and mysteries surrounding time reversal for global attractors of even the simplest of parabolic PDEs -- which diffusion, supposedly, governs ever so ``irreversibly''.

\textbf{Acknowledgment.}
This work has been supported, most generously, by the Deutsche Forschungsgemeinschaft, Collaborative Research Center 910 \emph{``Control of self-organizing nonlinear systems: Theoretical methods and concepts of application''} under project A4: \emph{``Spatio-temporal patterns: control, delays, and design.''}
We are grateful for the numerous inspirations, lively discussions, and excitingly active working atmosphere to which our speakers Eckehard Schöll and Sabine Klapp contributed so much.
We are also much indebted for enlightening discussions on meanders with very patient Piotr Zograf, and for the warm hospitality at the Mathematical Institute of Sankt Petersburg University.
Support by FCT/Portugal through projects UID/MAT/04459/2019 and UIDB/04459/2020 is also gratefully acknowledged.


\end{document}